\documentclass[11pt, reqno]{amsart}

\usepackage[OT2, T1]{fontenc}

\usepackage[usenames,dvipsnames]{color}

\usepackage{url}
\usepackage{amsmath}
\usepackage{array}
\usepackage{graphicx}
\usepackage{amssymb}
\usepackage{amsthm}
\definecolor{link}{RGB}{11,0,128}
\usepackage[colorlinks=true,citecolor=NavyBlue,linkcolor=Bittersweet,urlcolor=link]{hyperref}
\usepackage{hyperref}
\usepackage{paralist}
\usepackage{colonequals}		
\usepackage{color}
\usepackage{mathabx}		
\usepackage{amscd}
\usepackage[all,cmtip]{xy}	
\usepackage{sseq}			
\usepackage{verbatim}
\usepackage{parskip}
\usepackage{microtype}
\usepackage[in]{fullpage}
\usepackage{mathrsfs} 			
\usepackage{cleveref}
\usepackage{enumitem}
\usepackage{moreenum}			
\usepackage[alphabetic,lite]{amsrefs} 	
\usepackage{stmaryrd}	
\usepackage[foot]{amsaddr}
\usepackage{subfiles}
\usepackage{ragged2e}

\DeclareSymbolFont{cyrletters}{OT2}{wncyr}{m}{n}
\DeclareMathSymbol{\Sha}{\mathalpha}{cyrletters}{"58}

\newcommand{\bQ}{\mathbb{Q}}

\newcommand{\bZ}{\mathbb{Z}}

\newcommand{\cX}{\mathcal{X}}

\newcommand{\fa}{\mathfrak{a}}

\newcommand{\fm}{\mathfrak{m}}

\newcommand{\fp}{\mathfrak{p}}
\newcommand{\fq}{\mathfrak{q}}
\newcommand{\fr}{\mathfrak{r}}

\newcommand{\sH}{\mathscr{H}}
\newcommand{\sI}{\mathscr{I}}
\newcommand{\sJ}{\mathscr{J}}

\newcommand{\sL}{\mathscr{L}}
\newcommand{\sM}{\mathscr{M}}

\newcommand{\sO}{\mathscr{O}}

\newcommand{\ra}{\rightarrow}
\newcommand{\la}{\leftarrow}
\newcommand{\Ra}{\Rightarrow}

\newcommand{\xra}{\xrightarrow}

\newcommand{\hra}{\hookrightarrow}

\newcommand{\wt}{\widetilde}
\newcommand{\wh}{\widehat}

\newcommand{\ce}{\colonequals}
\newcommand{\ov}{\overline}
\newcommand{\un}{\underline}

\renewcommand{\b}{\textbf}
\newcommand{\surjects}{\twoheadrightarrow}

\newcommand{\tensor}{\otimes} 		
\newcommand{\isomto}{\overset{\sim}{\longrightarrow}}

\renewcommand{\i}{^{-1}}

\providecommand{\abs}[1]{\left\lvert#1\right\rvert}

\providecommand{\p}[1]{\left(#1\right)}
\providecommand{\up}[1]{{\upshape(}#1{\upshape)}}
\providecommand{\uref}[1]{{\upshape\ref{#1}}}
\providecommand{\uS}{{\upshape\S}}
\providecommand{\ucolon}{{\upshape:} }
\providecommand{\uscolon}{{\upshape;} }
\providecommand{\Sone}{\text{{\upshape($S_1$)}} }
\providecommand{\Stwo}{\text{{\upshape($S_2$)}} }
\providecommand{\Sn}{\text{{\upshape($S_n$)}} }
\providecommand{\Sonen}{\text{{\upshape($S_1$)}}}
\providecommand{\Stwon}{\text{{\upshape($S_2$)}}}
\providecommand{\Snn}{\text{{\upshape($S_n$)}}}

\providecommand{\f}[2]{\frac{#1}{#2}}
\DeclareMathOperator{\CM}{CM}		
\DeclareMathOperator{\Reg}{Reg}		
\DeclareMathOperator{\im}{Im}			
\DeclareMathOperator{\Spec}{Spec}		
\DeclareMathOperator{\Hom}{Hom}			
\DeclareMathOperator{\Ann}{Ann}			
\DeclareMathOperator{\depth}{depth}		
\DeclareMathOperator{\Supp}{Supp}		
\DeclareMathOperator{\Bl}{Bl}			
\newcommand{\ba}{\begin{aligned}}
\newcommand{\ea}{\end{aligned}}
\newcommand{\be}{\begin{equation}}
\newcommand{\ee}{\end{equation}}
\newcommand{\bpf}{\begin{proof}}
\newcommand{\epf}{\end{proof}}
\newcommand{\bthm}{\begin{thm}}
\newcommand{\ethm}{\end{thm}}
\newcommand{\bprop}{\begin{prop}}
\newcommand{\eprop}{\end{prop}}
\newcommand{\bcor}{\begin{cor}}
\newcommand{\ecor}{\end{cor}}
\newcommand{\brem}{\begin{rem}}
\newcommand{\erem}{\end{rem}}

\newcommand{\brems}{\begin{rems} \hfill \begin{enumerate}[label=\b{\thenumberingbase.},ref=\thenumberingbase]}
\newcommand{\remi}{\addtocounter{numberingbase}{1} \item}
\newcommand{\erems}{\end{enumerate} \end{rems}}

\newcommand{\begs}{\begin{egs} \hfill \begin{enumerate}[label=\b{\thenumberingbase.},ref=\thenumberingbase]}

\newcommand{\eegs}{\end{enumerate} \end{egs}}

\newcommand{\blem}{\begin{lemma}}
\newcommand{\elem}{\end{lemma}}

\newcommand{\bconj}{\begin{conj}}
\newcommand{\econj}{\end{conj}}
\newcommand{\bprob}{\begin{Problem}}
\newcommand{\eprob}{\end{Problem}}

\newcommand{\bq}{\begin{Q}}
\newcommand{\eq}{\end{Q}}
\newcommand{\benum}{\begin{enumerate}[label={{\upshape(\alph*)}}]}
\newcommand{\benuma}{\begin{enumerate}[label={{\upshape(\arabic*)}}]}
\newcommand{\benumr}{\begin{enumerate}[label={{\upshape(\roman*)}}]}
\newcommand{\eenum}{\end{enumerate}}
\newcommand{\bitem}{\begin{itemize}}
\newcommand{\eitem}{\end{itemize}}
\newcommand{\bc}{}
\newcommand{\bd}{\begin{defn}}
\newcommand{\ed}{\end{defn}}

\newcommand{\beg}{\begin{eg}}
\newcommand{\eeg}{\end{eg}}

\newcommand{\bcl}{\begin{claim}}
\newcommand{\ecl}{\end{claim}}
\newcommand{\lab}{\label}

\newcommand{\x}{\text}

\newcommand{\q}{\quad}
\newcommand{\qq}{\quad\quad}
\newcommand{\qqq}{\quad\quad\quad}

\newcommand{\tst}{\textstyle}

\newcommand*{\QED}{\hfill\ensuremath{\qed}}
\newcommand*{\QEDD}{\hfill\ensuremath{\qed\qed}}

\usepackage{aliascnt}
\newaliascnt{numberingbase}{subsection}

\theoremstyle{plain}
\newtheorem{thm}[numberingbase]{Theorem}
\Crefname{thm}{Theorem}{Theorems}

\Crefname{variant}{Variant}{Variants}

\Crefname{rethm}{Theorem}{Theorem}
\newtheorem{prop}[numberingbase]{Proposition}
\Crefname{prop}{Proposition}{Propositions} 
\newtheorem{Q}[numberingbase]{Question}
\Crefname{Q}{Question}{Questions}
\newtheorem{Problem}[subsection]{Problem}
\Crefname{Problem}{Problem}{Problems}
\newtheorem{conj}[numberingbase]{Conjecture}
\Crefname{conj}{Conjecture}{Conjectures}
\newtheorem{cor}[numberingbase]{Corollary}
\Crefname{cor}{Corollary}{Corollaries}
\newtheorem{lemma}[numberingbase]{Lemma}

\Crefname{subprop}{Proposition}{Propositions}

\Crefname{subcor}{Corollary}{Corollaries}

\Crefname{sublem}{Lemma}{Lemmas}

\theoremstyle{remark}
\newtheorem{claim}[equation]{Claim}
\Crefname{claim}{Claim}{Claims}

\Crefname{subrem}{Remark}{Remarks}

\theoremstyle{definition}
\newtheorem{defn}[numberingbase]{Definition}
\Crefname{defn}{Definition}{Definitions}

\Crefname{conv}{Convention}{Conventions}
\newtheorem{eg}[numberingbase]{Example}
\Crefname{eg}{Example}{Examples}
\newtheorem{rem}[numberingbase]{Remark}
\Crefname{rem}{Remark}{Remarks}
\newtheorem*{rems}{Remarks}
\newtheorem*{egs}{Examples}

\newtheoremstyle{subsection-tweak}
   {11pt}
   {3pt}%
   {}
   {}%
   {\bfseries}
   {}%
   {.5em}
   {\thmnumber{\@{#1}{}\@{#2}.}%
    \thmnote{~{\bfseries#3.}}}    
    
\theoremstyle{subsection-tweak}
\newtheorem{pp}[numberingbase]{}
\newcommand{\bpp}{\begin{pp}}
\newcommand{\epp}{\end{pp}}

\numberwithin{equation}{numberingbase}

\makeatletter
\def\@tocline#1#2#3#4#5#6#7{
    \begingroup 
    \@ifempty{#4}{%
    }{%
    }%

    \parindent\z@ \leftskip#3\relax \advance\leftskip\@tempdima\relax
    #5\hskip-\@tempdima
      \ifcase #1
       \or\or \hskip 2em \or \hskip 1em \else \hskip 3em \fi%
      #6\nobreak\relax
    \dotfill\hbox to\@pnumwidth{\@tocpagenum{#7}}\par
    \nobreak
    \endgroup
  }
 \def\l@section{\@tocline{1}{0pt}{1pc}{}{}}

\renewcommand{\tocsection}[3]{%
  \indentlabel{\@ifnotempty{#2}{\makebox[1.3em][l]{%
    \ignorespaces#1 \bfseries{#2}.\hfill}}}\bfseries{#3}
    \vspace{1.5pt}}

\renewcommand{\tocsubsection}[3]{%
  \indentlabel{\@ifnotempty{#2}{\hspace*{-0.5em}\makebox[2.1em][l]{%
    \ignorespaces#1#2.\hfill}}}#3
    \vspace{1.5pt}}

\makeatother

\begin{document}

\title{Macaulayfication of Noetherian schemes}

\author{K\k{e}stutis \v{C}esnavi\v{c}ius}
\address{CNRS, UMR 8628, Laboratoire de Math\'{e}matiques d'Orsay, Universit\'{e} Paris-Saclay, 91405 Orsay, France}
\email{kestutis@math.u-psud.fr}

\date{\today}


\begin{abstract} To reduce to resolving Cohen--Macaulay singularities, Faltings initiated the program of  ``Macaulayfying'' a given Noetherian scheme $X$. For a wide class of $X$, Kawasaki built the sought Cohen--Macaulay modifications, with a crucial drawback that his blowing ups did not preserve the locus $\CM(X) \subset X$ where $X$ is already Cohen--Macaulay. We extend Kawasaki's methods to show that every quasi-excellent, Noetherian scheme $X$ has a Cohen--Macaulay $\wt{X}$ with a proper map $\wt{X} \ra X$ that is an isomorphism over $\CM(X)$. This completes Faltings' program, reduces the conjectural resolution of singularities to the Cohen--Macaulay case, and implies that every proper, smooth scheme over a number field has a proper, flat, Cohen--Macaulay model over the ring of integers.
 \end{abstract}


\maketitle

\hypersetup{
    linktoc=page,     
}

\renewcommand*\contentsname{}
\q\\
\tableofcontents


 
 \section{Macaulayfication as a weak form of resolution of singularities}

The resolution of singularities, in Grothendieck's formulation, predicts the following.

\bconj \lab{res-sing-conj}
For a quasi-excellent,\footnote{We recall from \cite{EGAIV2}*{7.8.2, 7.8.5} and \cite{ILO14}*{I.2.10, I.7.1} that a scheme $X$ is \emph{quasi-excellent} if it is locally Noetherian, for each $x \in X$ the geometric fibers of the map 
\[
\Spec(\wh{\sO}_{X,\, x}) \ra \Spec(\sO_{X,\, x})
\]
are regular (the fibers of this map are the \emph{formal fibers} of $X$ at $x$), and every integral $X$-scheme $X'$ that is finite over some affine open of $X$ has a nonempty regular open subscheme $U' \subset X'$. For such an $X$, the subset $\Reg(X)$ of the points $x \in X$ for which the local ring $\sO_{X,\, x}$ is regular is open. A quasi-excellent $X$ is \emph{excellent} if it is universally catenary (see \S\ref{catenary}).
} reduced, Noetherian scheme $X$, there are a regular scheme $\wt{X}$ and a proper morphism $\pi \colon \wt{X} \ra X$ that is an isomorphism over the regular locus $\Reg(X)$ of $X$.
\econj

The conjecture is interesting even without requiring $\pi|_{\pi\i(\Reg(X))}$ to be an isomorphism, but this condition is natural: for instance, if one wishes to find a proper and regular integral model of a proper and smooth scheme $Y$ over a number field, one does not want to have to modify $Y$. 

When $X$ is a $\bQ$-scheme \Cref{res-sing-conj} is known by Temkin's \cite{Tem08}*{Thm.~1.1}, which builds on Hironaka's work \cite{Hir64} in a slightly more restrictive setting (such a generalization of \emph{op.~cit.}~was claimed in \cite{EGAIV2}*{7.9.6}).  When $X$ is of positive or mixed characteristic, the conjecture is open apart from low dimensional cases and variants, such as alterations \cite{dJ96}. In general, one may reduce \Cref{res-sing-conj} to the case when $X$ is excellent and locally equidimensional, see \Cref{RS-exc-red}.

It is natural to approach resolution of singularities by trying to gradually improve $X$: for instance, as initiated by  Faltings \cite{Fal78}, one may seek a weaker version of \Cref{res-sing-conj} in which the regularity of $\wt{X}$ is weakened to Cohen--Macaulayness. It seems difficult to achieve this with the Hironaka style methods based on blowing up regular centers: for example, to Macaulayfy a Buchsbaum singularity in any characteristic one blows up the possibly nonreduced ideal generated by a system of parameters.

The main purpose of the present article is to establish the variant proposed by Faltings, and hence to reduce the resolution of singularities to the Cohen--Macaulay case (for previous work on this, see \S\ref{prev-work}). In fact, for this, quasi-excellence, whose definition is modeled on regularity, is unnaturally restrictive. We will replace it by the following weaker CM-quasi-excellence requirement.

\bd \lab{CM-exc}
A scheme $X$ is \emph{CM-quasi-excellent} if it is locally Noetherian and such that
\benuma
\item \lab{CME-1}
the formal fibers of the local rings of $X$ are Cohen--Macaulay\uscolon

\item \lab{CME-2}
every integral, closed subscheme $X' \subset X$ has a nonempty, Cohen--Macaulay open subscheme\uscolon
\eenum
a CM-quasi-excellent $X$ is \emph{CM-excellent} if, in addition, it is universally catenary.
\ed

CM-quasi-excellence (even \ref{CME-2} alone) implies the openness of the Cohen--Macaulay locus $\CM(X)$ (see \S\ref{open-Sn}).

\beg \lab{CM-exc-eg}
Every quasi-excellent scheme is CM-quasi-excellent, and similarly for excellence. In contrast to the existence of nonexcellent discrete valuation rings caused by the possible nonseparability of the field extension obtained by completing the fraction field, every \emph{Dedekind scheme}, that is, a Noetherian normal scheme of dimension $\le 1$, is CM-excellent. 
\eeg

\brem
\lab{quasi-CM}
In fact, every Cohen--Macaulay scheme is CM-excellent. More generally, if a locally Noetherian scheme $X$ is \emph{quasi-Cohen--Macaulay} in the sense that it has a coherent Cohen--Macaulay $\sO_X$-module $\sM$ with $\Supp(\sM) = X$, then any closed subscheme of $X$ is CM-excellent: \cite{EGAIV2}*{6.3.8 and 6.11.9~(ii)} ensure \ref{CME-1} and \ref{CME-2}, whereas the universal catenarity follows from the proof of \cite{EGAIV2}*{6.3.7} that continues to work in the quasi-Cohen--Macaulay case.
\erem

\brem\lab{CM-exc-stab}
 If $X$ is CM-quasi-excellent, then so is every localization of a locally finite type $X$-scheme thanks to \cite{EGAIV2}*{7.4.4 (with 7.3.8)} for \ref{CME-1} and \cite{EGAIV2}*{6.11.9~(i)} for \ref{CME-2}. 
\erem

Our main result is the promised Cohen--Macaulay version of \Cref{res-sing-conj}:

\bthm[\Cref{main-thm-pf}, \Cref{rem-proj}] \lab{main-thm}
For every CM-quasi-excellent, Noetherian scheme $X$, there are a Cohen--Macaulay scheme $\wt{X}$ and a birational, projective morphism 
\[
\pi\colon \wt{X} \ra X
\]
that is an isomorphism over the Cohen--Macaulay locus 
$\CM(X) \subset X$. 
\ethm

\brem
For any fixed closed subscheme $Z \subset X$ whose support is $X\setminus \CM(X)$, we may arrange that, in addition, the (scheme-theoretic) preimage of $Z$ in $\wt{X}$ be a divisor (with support $\wt{X} \setminus \CM(X)$). Indeed, since $\pi$ is birational and $\wt{X}$ has no embedded associated primes, the preimage of any divisor in $X$ is a divisor in $\wt{X}$, so it suffices to apply \Cref{main-thm} with $\Bl_Z(X)$ in place of $X$.
\erem


\Cref{main-thm} is sharp in the sense that if each integral, closed subscheme of a locally Noetherian scheme $X$ admits a Macaulayfication, then $X$ must be CM-quasi-excellent, see \Cref{converse}.

Since CM-quasi-excellence is weaker than quasi-excellence, \Cref{main-thm} is new even for $\bQ$-schemes. Nevertheless, it is of most interest in positive and mixed characteristic. For instance, it implies the existence of proper, flat, Cohen--Macaulay integral models over rings of integers of number fields:

\bcor \lab{int-models}
For every integral Dedekind scheme $S$ with the function field $K$ and every proper, Cohen--Macaulay $K$-scheme $X$, there is a proper, flat, Cohen--Macaulay $S$-scheme $\cX$ with $\cX_K \simeq X$. If $X$ is projective over $K$, then one may choose $\cX$ to be projective over $S$.
\ecor

Since $S$ is CM-excellent (see \Cref{CM-exc-eg}), \Cref{int-models} follows by applying \Cref{main-thm}  to the schematic image $\cX_0$ of $X$ in a Nagata compactification over $S$ (to use the Nagata compactification \cite{Del10}*{1.6}, one first spreads out $X$): indeed, since $\cX_0$ has $X$ as the $K$-fiber and is $S$-flat, 
so is the resulting $\cX$; 
for the projective aspect, one instead forms the schematic image in a projective space. 

\brem
If $X$ in \Cref{int-models} comes equipped with a finite family of coherent Cohen--Macaulay modules (for instance, vector bundles), then, by using the finer version of \Cref{main-thm} stated in \Cref{main-thm-pf}, one may arrange that, in addition, they extend to Cohen--Macaulay modules on $\cX$.
\erem

Another basic consequence of \Cref{main-thm} (and Remark \ref{quasi-CM}) is the following principalization result:

\bcor \lab{principalize}
For every Noetherian, Cohen--Macaulay scheme $X$ and every closed subscheme $Z \subset X$, there are a Cohen--Macaulay scheme $\wt{X}$ and a projective morphism $\wt{X} \ra X$ such that the \up{scheme-theoretic} preimage of $Z$ in $\wt{X}$ is a divisor and $\wt{X} \ra X$ is an isomorphism over the maximal open subscheme $U \subset X$ on which $Z$ is already a divisor.
\ecor

To obtain \Cref{principalize}, one first notes that both $X$ and $\Bl_Z(X)$ are CM-excellent and locally equidimensional (see \Cref{fed-stab}~\ref{FS-b}), the map $\Bl_Z(X) \ra X$ is an isomorphism over $U$, and the preimage of $Z$ in $\Bl_Z(X)$ is a divisor. One then applies \Cref{main-thm} (in its more precise form \ref{main-thm-pf}) to $\Bl_{Z}(X)$ to obtain a Macaulayfying blowing up $\wt{X} \ra \Bl_{Z}(X)$ that is an isomorphism over $U$. The preimage of $Z$ in $\wt{X}$ is locally principal and, as is checked away from the exceptional divisor of $\wt{X} \ra \Bl_{Z}(X)$, is a divisor.

By combining the Nagata compactification with \Cref{main-thm} and \Cref{principalize}, we obtain the following consequence that concerns the existence of Cohen--Macaulay compactifications.

\bcor\lab{compactify}
For every CM-quasi-excellent Noetherian scheme $S$ and every finite type, separated $S$-scheme $X$ that is Cohen--Macaulay, there is an open $S$-immersion $X \hra \ov{X}$ into a proper $S$-scheme $\ov{X}$ that is Cohen--Macaulay such that $\ov{X} \setminus X$ is a \up{possibly nonreduced} divisor in $X$. 
\ecor

\bpp[Previous work on Macaulayfication] \lab{prev-work}
The key novelty of \Cref{main-thm} is that its Macaulayfication map preserves the Cohen--Macaulay locus of $X$ as is crucial for the corollaries above. Previously this has only been achieved in cases when the non-Cohen--Macaulay locus 
\[
X \setminus \CM(X)
\]
is a disjoint union of points, see \cite{Fal78}*{Satz 2} and \cite{Bro83a}*{\S6.C}, a situation intimately related to the study of generalized Cohen--Macaulay rings (such as Buchsbaum rings) and of their Macaulayfications, see \cite{Got82}*{1.1} and \cite{Sche83}*{4.2}. Without this isomorphy condition over $\CM(X)$ and under suitable additional assumptions, birational Macaulayfications were constructed in \cite{Fal78}*{Satz 3} and \cite{Bro83b}*{1.1, 1.2} when $\dim(X \setminus \CM(X)) \le 1$, in \cite{Kaw98}*{1.1} when $\dim(X \setminus \CM(X)) \le 2$, and in \cite{Kaw00}*{1.1} with complements in \cite{Kaw02}*{1.1} and \cite{Kaw08}*{1.2} when $\dim(X \setminus \CM(X))$ is arbitrary; another method, based on the idea of de Jong's alterations, was exhibited in \cite{Hon04}*{\S4.1} (see also \cite{Hei14}). Macaulayfication of coherent modules was explored in \cite{Mor99}.
\epp

\bpp[The inductive method]
Our technique extends that of Kawasaki used in \cite{Kaw00}, which in turn builds on the one of Faltings used in \cite{Fal78}. After an initial reduction to locally equidimensional schemes based on an inductive construction of an ``\Stwon-ificaiton,'' 
we use Noetherian induction to reduce to Macaulayfying a projective $X$-scheme $Y \ra X$ in the case when $X$ is the spectrum of a complete, Noetherian, local ring and $x \in X$ is the closed point. By induction, $Y \setminus Y_x$ is already Cohen--Macaulay, and we need to find a Macaulayfying blowing up with the center contained in a thickening of $Y_x$. For this, it is key to allow the center to meet $\CM(Y) \cap Y_x$: for instance, to resolve $Y_x$ itself, we would choose a Kawasaki-style center constructed from well-chosen hypersurfaces that cut out a Cohen--Macaulay global complete intersection in $Y_x$. Since instead we need to resolve $Y$, we first make a preliminary blowing up to make $Y_x \subset Y$ into a divisor; we then choose hypersurfaces on $Y$ whose restrictions to $Y_x$ are like in Kawasaki's method. The key trick is to complement these with a large power of the ideal $\sI_{Y_x}$ considered as an additional hypersurface, and then to build a Kawasaki-style center from this larger collection. This is legitimate because $\sI_{Y_x}$ is locally principal and Cohen--Macaulayness is also local. By regarding the power of $\sI_{Y_x}$ as the ``first'' hypersurface in the collection, we can keep the constructed center disjoint from $Y \setminus Y_x$.
\epp

\bpp[Notation and conventions] \lab{conv}
For a coherent module $\sM$ on a locally Noetherian scheme $X$, its \emph{support} is the closed subscheme $\Supp(\sM) \subset X$ cut out by the annihilator ideal $\Ann_{\sO_X}(\sM) \subset \sO_X$ (the latter is coherent because $\sM$ is finitely generated). For $n \in \bZ$, such an $\sM$  is \emph{(S$_n$)} if
\[
\depth_{\sO_x}(\sM_x) \ge \min(n, \dim(\Supp(\sM_x))) \q \x{for all} \q x \in X
\]
(we recall from \cite{EGAIV1}*{0.14.1.2} that $\dim(\emptyset) = -\infty$). For instance, $\sM$ is \Sone if and only if it has no embedded associated primes (see \cite{EGAIV2}*{5.7.5}). Moreover, $\sM$ is \emph{Cohen--Macaulay} if it is \Sn for every $n$, that is, if 
\[
\depth_{\sO_x}(\sM_x) = \dim(\Supp(\sM_x)) \q \x{for all}\q x \in X.
\]
A scheme $X$ is \emph{(S$_n$)} or \emph{Cohen--Macaulay} if it is locally Noetherian and $\sO_X$ has the respective property as an $\sO_X$-module. We let $\abs{X}$ denote the underlying topological space of a scheme $X$.

For a scheme $X$, the \emph{height} of an $x \in X$ is the dimension of the stalk $\dim(\sO_{X,\, x})$; the \emph{coheight} of $x$ is the dimension of the closure $\ov{\{ x\}} \subset X $. We denote the subset of points of height $i$ by $X^{(i)}$. The \emph{codimension} of a closed subscheme $Y \subset X$ is the infimum of the heights of points of $X$ that lie on $Y$ (compare with \cite{EGAIV2}*{5.1.3}). We denote the $X$-scheme obtained by blowing up $Y$ by $\Bl_Y(X)$; for the quasi-coherent ideal $\sI \subset \sO_X$ that cuts out $Y$, we also write $\Bl_{\sI}(X)$ for $\Bl_Y(X)$. We freely use that $\Bl_Y(X)$ has a universal property (see \cite{SP}*{\href{https://stacks.math.columbia.edu/tag/0806}{0806}}), commutes with flat base change in $X$ (see \cite{SP}*{\href{https://stacks.math.columbia.edu/tag/0805}{0805}}), and is projective over $X$ when $\sI$ is of finite type (see \cite{EGAII}*{5.5.2}). 

For a module $M$ over a commutative ring $R$ and an ideal $\fr \subset R$, we write $M\langle\fr\rangle$ for the $\fr$-torsion submodule. We write $M\langle\fr^\infty\rangle$ for $\bigcup_{n > 0} M\langle\fr^n\rangle$ and simply $M\langle r\rangle$, etc.~when $\fr = (r)$ is principal. For a submodule $M' \subset M$, we use the colon notation $M' :_M \fr$ to denote the preimage of $(M/M')\langle\fr\rangle$ in $M$. If $R$ is local and Noetherian, we let $\wh{R}$ denote its completion with respect to the maximal ideal.
\epp

\subsection*{Acknowledgements}
I thank Cuong Trung Doan, David Hansen, Luc Illusie, Takesi Kawasaki, David Rydh, and Olivier Wittenberg for helpful conversations and correspondence. I thank the referees for helpful suggestions. I thank the CNRS and the Universit\'{e} Paris-Sud for their support during my work on this article. This project received funding from the European Research Council (ERC) under the European Union's Horizon 2020 research and innovation programme (grant agreement No.~851146).


\section{(S$_2$)-ification of coherent modules} \lab{S2ify-section}

The principal goal of this section is to reduce our search of a Macaulayfication to the case of locally equidimensional and CM-excellent (in particular, universally catenary) schemes, a case to which we will extend Kawasaki's approach in \S\ref{local-case}. This is a natural, even if not an entirely canonical, initial reduction because every Cohen--Macaulay (or even quasi-Cohen--Macaulay) scheme  satisfies these conditions (see \cite{EGAIV1}*{0.16.5.4} and Remark \ref{quasi-CM}).  For the resolution of singularities conjecture \ref{res-sing-conj}, the corresponding reduction is simpler, so our first goal is to explain it in \Cref{RS-exc-red}.

\bpp[Catenarity] \lab{catenary}
We recall from \cite{EGAIV2}*{5.6.3 (ii)} that a scheme $X$ is \emph{universally catenary} if it is locally Noetherian and every scheme $X'$ that is locally of finite type over $X$ is catenary in the sense that any two saturated chains of specializations of points of $X'$ with the same endpoints have the same length (any such chain is contained in every affine open of $X'$ that contains the endpoint of larger height). 
We already mentioned in Remark \ref{quasi-CM} that every closed subscheme of a quasi-Cohen--Macaulay scheme, is universally catenary. This implies, in particular, that every complete, Noetherian, local ring is universally catenary.
\epp

\bpp[Equidimensionality] \lab{formal-eq}
We say that a scheme $X$ is \emph{locally equidimensional} if it is locally Noetherian and each of its local rings $\sO_{X, x}$ is equidimensional in the sense that all the irreducible components of $\Spec(\sO_{X, x})$ have the same dimension. Moreover, similarly to \cite{EGAIV2}*{7.1.1}, we say that a scheme $X$ is \emph{formally equidimensional} (or \emph{quasi-unmixed} in other terminology) if it is locally Noetherian and the completion $\wh{\sO}_{X,\, x}$ is equidimensional for every $x \in X$. We recall from \cite{HIO88}*{18.17} that, by a result of Ratliff, a locally Noetherian scheme $X$ is formally equidimensional if and only if it is locally equidimensional and universally catenary.
\epp

The principal advantage of formal equidimensionality for the purpose of improving singularities is its pleasant interaction with blowing up. We now review this critical ingredient to our arguments.

\blem \lab{fed-stab}
Let $X$ be a formally equidimensional scheme and let $\sI \subset \sO_X$ be a coherent ideal.
\benum
\item \lab{FS-a}
If $\sI$ is locally principal of height $> 0$, then $\un{\Spec}(\sO_X/\sI)$ is formally equidimensional.

\item \lab{FS-b}
The blowing up $\Bl_{\sI}(X)$ is formally equidimensional.
\eenum
\elem

\bpf
Part \ref{FS-a} is a special case of \cite{HIO88}*{18.20}. Alternatively: the maximal saturated chains of primes  of $\wh{\sO}_{X,\, x}$ all have the same length, the completion of the local ring at $x \in \un{\Spec}(\sO_X/\sI)$ is $\wh{\sO}_{X,\, x}/\sI_x \wh{\sO}_{X,\, x}$, and, by \cite{EGAIV2}*{2.3.4}, the principal ideal $\sI_x \wh{\sO}_{X,\, x} \subset \wh{\sO}_{X,\, x}$ either vanishes or any minimal prime containing it has height $1$, so the maximal saturated chains of primes of $\wh{\sO}_{X,\, x}/\sI_x \wh{\sO}_{X,\, x}$ also all have the same length.  Part \ref{FS-b} follows from \cite{HIO88}*{18.26}. \epf

The following standard lemma will be useful for us on several occasions.

\blem \lab{Sn-descend}
For a Noetherian, local ring $R$ with \Sn formal fibers and an $R$-scheme $X$ such that both $X$ and $X_{\wh{R}}$ are locally Noetherian, a coherent $\sO_X$-module $\sM$ is \Sn if and only if\,\footnote{Even though the `if' direction does not require the assumption on the formal fibers (see \cite{EGAIV2}*{6.4.1~(i)}), the `only if' does even when $n = 1$: indeed, \cite{FR70}*{3.3} 
(which answered a question raised in \cite{EGAIV2}*{6.4.3}) exhibited a $2$-dimensional Noetherian local ring that is \Sone (even a domain) whose completion is not \Sonen.
} so is $\sM|_{X_{\wh{R}}}$. 
\elem

\bpf
Each fiber of the flat map $X_{\wh{R}} \ra X$ is the base changes of a fiber of $\Spec(\wh{R}) \ra \Spec(R)$ to a possibly larger field, and hence is \Snn. Thus, the claim is a special case of \cite{EGAIV2}*{6.4.2}.
\epf

Thanks to the result of Ratliff reviewed in \S\ref{formal-eq}, the following lemma is a useful source of catenarity.

\blem\lab{miracle-catenary}
An  \Stwo scheme whose local rings have \Stwo formal fibers is formally equidimensional. 

\elem

\bpf
By \Cref{Sn-descend}, the completions of the local rings are \Stwon. By \S\ref{catenary}, they are also catenary. Thus, \cite{EGAIV2}*{5.10.9} implies that these completions are equidimensional.
\epf

\brem
The assumption on the formal fibers cannot be dropped in \Cref{miracle-catenary}: indeed, Noetherian, noncatenary, normal, local domains exist by \cite{Ogo80}*{Appendix} or \cite{Hei82}. 
\erem

We are ready for the reduction of \Cref{res-sing-conj} to the excellent and locally equidimensional case.

\bprop \lab{RS-exc-red}
For a quasi-excellent, reduced scheme $X$, the normalization morphism $\wt{X} \ra X$ is a finite map that is an isomorphism over $\Reg(X)$, and $\wt{X}$ is excellent and formally equidimensional.
\eprop

\bpf
The $X$-finiteness of $\wt{X}$ is a minor improvement to \cite{EGAIV2}*{7.8.6~(ii)}, which was written for excellent $X$ (see \cite{ILO14}*{I, \S6} for a stronger such improvement). To obtain it, we first note that the Nagata criterion \cite{EGAIV2}*{7.7.3} ensures that the coordinate rings of the affine opens of $X$ are universally Japanese and then apply the definition of the normalization \cite{EGAII}*{6.3.4, 6.3.8}. Since $\wt{X}$ inherits quasi-excellence, \Cref{miracle-catenary} implies its formal equidimensionality, and so excellence. 
\epf

For the analogous reduction of the Macaulayfication problem to the CM-excellent case, the principal complication is the absence of a canonical ``\Stwon-ification'' morphism that would replace the normalization (for instance, \cite{Bro86}*{3.9--3.11} confirms such absence). Indeed, even for \Sone schemes, the naive approach of pushing forward the structure sheaf from the open \Stwo locus does not work because such pushforward may fail to be coherent: this happens, for instance, in the case of a $2$-dimensional, Noetherian, local ring that has an irreducible component of dimension $1$. We will build a noncanonical \Stwon-ification in \Cref{S2ify-sch} after the following preparations.

\bpp[Openness of the (S$_n$)-loci] \lab{open-Sn}
For a locally Noetherian scheme $X$ and a coherent $\sO_X$-module $\sM$, the subset 
\[
U_{(\x{S}_1)}(\sM) \subset X
\]
of points at which the stalk of $\sM$ is (S$_1$) is open: indeed, the restriction of $U_{(\x{S}_1)}(\sM)$ to any affine open of $X$ is the complement of the union of the closed subschemes cut out by the embedded associated primes of $\sM$ (see \cite{EGAIV2}*{6.11.7 (i) and its proof}). In contrast, the subset 
\[
U_{(\x{S}_n)}(\sM) \subset X
\]
of points at which the stalk of $\sM$ is (S$_n$) need not be open for $n > 1$, see \cite{FR70}*{3.5} (which answered the question raised in \cite{EGAIV2}*{6.11.9 (ii)}). However, \cite{EGAIV2}*{6.11.6} ensures that $U_{(\x{S}_n)}(\sM)$ is open in $X$ for every $\sM$ if every integral, closed subscheme $X' \subset X$ has a nonempty open subscheme that is {\upshape(S$_n$)}, for instance, if $X$ is CM-quasi-excellent. Similarly, \cite{EGAIV2}*{6.11.8} ensures that the subset 
\[
\CM(\sM) \subset X
\]
of points of $X$ at which the stalk of $\sM$ is Cohen--Macaulay is open for every coherent $\sO_X$-module $\sM$ if every $X'$ as above contains a nonempty open subscheme that is Cohen--Macaulay, for instance, if $X$ is CM-quasi-excellent.

For brevity, we often write $U_{\Sn}(X)$ and $\CM(X)$ in place of $U_{(\x{S}_n)}(\sO_X)$ and $\CM(\sO_X)$, respectively.
\epp

Bearing the openness of $U_{\Sone}(\sM)$ in mind, one may easily build an \Sonen-ification of $\sM$ as follows.

\bthm \lab{S1-ify}
For a coherent module $\sM$ on locally Noetherian scheme $X$,
\[
\tst \sM' \ce \im\p{\sM \ra  j_*(\sM|_{U_{\Sone}(\sM)})}, \quad \text{where} \quad j \colon U_{\Sone}(\sM) \hra X
\]
is the indicated open immersion, is an {\upshape(S$_1$)}-ification of $\sM$\ucolon it is a coherent $\sO_X$-module that is {\upshape(S$_1$)} and agrees with $\sM$ on $U_{\Sone}(\sM)$.
\ethm

\bpf
By construction, $\sM'$ is coherent, agrees with $\sM$ on $U_{\Sone}(\sM)$, and has no nonzero local sections that vanish on $U_{\Sone}(\sM)$. Thus, the supports of $\sM$ and $\sM'$ agree topologically (both are equal to the closure of the generic points of $\Supp(\sM)$) and $\sM'$ has no embedded associated primes, that is, is  \Sone (see \S\ref{conv}).
\epf

\bpp[\Snn-quasi-excellence] \lab{Sn-qe}
For $n \in \bZ$, a scheme $X$ is \emph{(S$_n$)-quasi-excellent} if it is locally Noetherian and such that

\benuma
\item \lab{Sn-qe-1}
the formal fibers of the local rings of $X$ are \Snn\uscolon

\item \lab{Sn-qe-2}
every integral, closed subscheme $X' \subset X$ has a nonempty, \Sn open subscheme;
\eenum
an \Snn-quasi-excellent $X$ is \emph{(S$_n$)-excellent} if, in addition, it is universally catenary.

By \S\ref{open-Sn}, condition \ref{Sn-qe-2} implies that the \Sn locus $U_{\x{\Snn}}(\sM) \subset X$ is open for every coherent $\sO_X$-module $\sM$.  Evidently, a CM-quasi-excellent (resp.,~CM-excellent) scheme is \Snn-quasi-excellent (resp.,~\Snn-excellent) for every $n$. The references used in Remark \ref{CM-exc-stab} imply that \Snn-quasi-excellence is stable under localization and ascends along morphisms that are locally of finite type.
\epp

We will build \Stwon-ifications of \Stwon-quasi-excellent schemes $X$ in \Cref{S2-ify}. For this, we begin with the following auxiliary lemma that treats the technically simpler case of \Stwon-excellent $X$.

\blem \lab{S2-ify-lem}
For an open immersion $j \colon U \hra X$ of Noetherian, \Stwon-excellent schemes, any coherent $\sO_U$-module $\sM$ that is \Sone \up{resp.,~\Stwon} extends to a coherent $\sO_X$-submodule 
\be \lab{S2L-eq}
\sM' \subset j_*(\sM)
\ee
that is \Sone \up{resp.,~a finite direct sum of \Stwo modules}\uscolon if $\sM$ underlies a commutative $\sO_U$-algebra, then  $\sM'$ may be chosen to be an algebra extension of $\sM$ \up{with the direct sum that of $\sO_X$-algebras}. Moreover, if $\sM$ is \Stwo and for each $x \in \ov{\Supp(\sM)} \setminus \Supp(\sM)$ of height $\ge 2$ in $\ov{\Supp(\sM)}$  the punctured spectrum of $\sO_{\ov{\Supp(\sM)},\, x}$ has no isolated points, then $\sM'$ may be chosen to itself be \Stwon.
\elem

\bpf
By replacing $X$ by the schematic image of $\Supp(\sM)$, we may assume that 
\[
\Supp(\sM) = U
\]
and $U$ is dense in $X$ (the formation of the schematic image commutes with localization by \cite{EGAI}*{9.5.8}). Noetherian induction on $X\setminus U$ and spreading out allow us to localize at a generic point of $X\setminus U$ to assume that $X$ is local, $U$ is the (nonempty) complement of the closed point, and, by passing to direct summands in the \Stwo case if needed, that $\sM$ is \Sone (resp.,~\Stwon, not merely a direct sum). 

In the case $\dim X = 1$, we use \cite{EGAI}*{9.4.7} to find a coherent $\sO_X$-submodule 
\[
\sM' \subset j_*(\sM)
\]
extending $\sM$; by construction, $\sM'$ has no embedded associated primes, so is \Sonen, and hence also \Stwon. If $\dim X = 1$ and $\sM$ is a commutative $\sO_U$-algebra, we proceed differently: we use Zariski's main theorem \cite{EGAIV4}*{18.12.13} and the fact that $\dim(U) = 0$ to extend $\sM$ to a commutative $\sO_X$-algebra $\sM''$ that is coherent as an $\sO_X$-module; we then let $\sM'$ be the image of the map 
\[
\sM'' \ra j_*(\sM)
\]
to obtain a desired commutative $\sO_X$-algebra extension that is \Sonen, so also \Stwon.

In the remaining case $\dim X \ge 2$, we write 
\[
\tst U \cong U_0 \bigsqcup U_{\ge 1}
\]
where $U_0$ (resp.~$U_{\ge 1}$) is the union of the isolated points (resp.,~of the irreducible components of dimension $\ge 1$) of $U$. The module $\sM$ decomposes accordingly: 
\[
\tst \sM \cong \sM|_{U_0} \bigoplus \sM|_{U_{\ge 1}}.
\]
The schematic image of $U_0$ in $X$ is $1$-dimensional and local, so the settled $\dim X = 1$ case supplies a desired extension of $\sM|_{U_0}$. We may therefore assume that 
\[
\Supp(\sM) = U_{\ge 1}
\]
(in the \Stwo case, this step is the reason for a direct sum in the statement; in the \Sone case, no direct sum issue arises because \eqref{S2L-eq} alone implies that $\sM'$ has no embedded associated primes and hence is \Sonen). For such $\sM$, we set 
\[
\sM' \ce j_*(\sM)
\]
and aim to show that the $\sO_X$-module $\sM'$ is coherent. It will then follow from \cite{EGAIV2}*{5.10.5} that $\sM'$ is \Stwon.

For the coherence of $j_*(\sM)$, we will use Koll\'{a}r's criterion \cite{SP}*{\href{http://stacks.math.columbia.edu/tag/0BK3}{0BK3}} (or \cite{Kol17}*{Thm.~2}) that gives a necessary and sufficient condition (for similar earlier results, see \cite{EGAIV2}*{5.11.4, 7.2.2} and \cite{SGA2new}*{VIII,~2.3}). We need to check that for every associated prime $u \in U$ of $\sM$ and the closed point $x$ of $X$, 
the coheights of the associated primes of $\wh{\sO}_{\ov{\{u \}},\, x}$ are all $\ge 2$. Since $\sO_{\ov{\{u \}},\, x}$ inherits \Stwon-quasi-excellence from $X$ (see \S\ref{Sn-qe}) and is a domain, by \Cref{Sn-descend}, its completion $\wh{\sO}_{\ov{\{u \}},\, x}$ is \Sonen, and hence has no embedded associated primes. Thus, since, by a result of Ratliff \cite{SP}*{\href{https://stacks.math.columbia.edu/tag/0AW6}{0AW6}}, the universal catenarity of $X$ implies that $\wh{\sO}_{\ov{\{u \}},\, x}$ is equidimensional, the coheights in question are all equal to $\dim(\sO_{\ov{\{u \}},\, x})$. Since $\sM$ is \Sone and $U_{\ge 1}$ is its support, $u$ is a generic point of $U_{\ge 1}$, so that 
\[
\dim(\sO_{\ov{\{u \}},\, x}) \ge 2. \qedhere
\]
\epf


\brem
The \Stwon-excellence of $X$ in \Cref{S2-ify-lem} may be weakened to the combination of \S\ref{Sn-qe}~\ref{Sn-qe-2} with $n = 2$ and strict formal catenarity.\footnote{Building on \cite{EGAIV2}*{7.2.1, 7.2.6}, we call a scheme $X$ \emph{strictly formally catenary} if it is locally Noetherian and every integral, closed subscheme of $X' \subset X$ is \emph{strictly formally equidimensional} in the sense that for every $x \in X'$ the completed local ring $\wh{\sO}_{X',\, x}$ is equidimensional and has no embedded associated primes.} Indeed, \S\ref{Sn-qe}~\ref{Sn-qe-2} ensures the openness of $U_{\x{\Stwon}}(\sM) \subset X$, whereas strict formal catenarity ensures that $\wh{\sO}_{\ov{\{u \}},\, x}$ is \Sone and equidimensional.
\erem


\bthm \lab{S2-ify}
For an open immersion $j\colon U \hra X$ of Noetherian, \Stwon-quasi-excellent schemes, an \Stwon, finite $U$-scheme $\wt{U}$, and finitely many \Stwon, coherent $\sO_{\wt{U}}$-modules $\sM_i$ with $\abs{\Supp(\sM_i)} = |\wt{U}|$, there are an \Stwon, finite $X$-scheme $\wt{X}$ and \Stwon, coherent $\sO_{\wt{X}}$-submodules $\sM'_i \subset (j|_{\wt{X}})_*(\sM_i)$ such~that 
\[
\wt{X}|_U \cong \wt{U} \q\x{and}\q \sM'_i|_{\wt{U}} = \sM_i
\]
and, in addition, $\abs{\Supp(\sM'_i)} = |\wt{X}|$ and  $\sO_{\wt{X}} \hra j_*(\sO_{\wt{U}})$, so that $\sO_{\wt{X}}$ is \Sone as an $\sO_X$-module. 
\ethm

\bpf
The inclusion $\sO_{\wt{X}} \hra j_*(\sO_{\wt{U}})$ shows that the associated primes of the $\sO_X$-module $\sO_{\wt{X}}$ coincide with those of the $\sO_U$-module $\sO_{\wt{X}}|_U \cong \sO_{\wt{U}}$, so our $\sO_{\wt{X}}$ will necessarily be \Sone as an $\sO_X$-module.

By replacing $X$ by the schematic image of $\wt{U}$, we may assume that $\Supp_{\sO_U}(\sO_{\wt{U}}) = U$ and $U$ is dense in $X$. As in the proof of \Cref{S2-ify-lem}, we may then localize at a generic point of $X\setminus U$ to assume that $X$ is local (necessarily of dimension $> 0$) and $U$ is the complement of the closed point. 

We then combine formal patching \cite{FR70}*{4.2} with the commutativity of $j_*(-)$ with flat base change to assume that $X$ is complete (to descend the coherence and the \Stwo property from the completion we use \cite{SGA1new}*{VIII, 1.10} and \Cref{Sn-descend}). In particular, since $X$ is now \Stwon-excellent, \Cref{S2-ify-lem} applied to the coordinate algebras of the connected components of $\wt{U}$ supplies a desired $\wt{X}$ such that its connected components correspond bijectively to those of $\wt{U}$ via pullback. Another application of \Cref{S2-ify-lem}, this time over $\wt{X}$, then supplies the desired \Stwo extensions $\sM'_i$ of the $\sM_i$. 
\epf

\bcor \lab{S2ify-sch}
Every Noetherian, \Stwon-quasi-excellent scheme $X$ has an \Stwon-ification\ucolon there is a finite birational map $\wt{X} \xra{\pi} X$ that is an isomorphism over $U_{\Stwo}(X)$ such that $\wt{X}$ is \Stwon, locally equidimensional, \Stwon-excellent, and has $U_{\Stwo}(X)$ as a dense open. 
\ecor

\bpf
\Cref{S2-ify} applied to $U = \wt{U} = U_{\Stwon}(X)$ supplies a candidate $\wt{X}$ that is \Stwo and for which $\sO_{\wt{X}}$ is \Sone as an $\sO_X$-module. Since $\wt{X}$ inherits \Stwon-quasi-excellence, \Cref{miracle-catenary} implies that it is formally equidimensional and \Stwon-excellent. Due to the \Sone property of $\sO_{\wt{X}}$, the generic points of $\wt{X}$ lie over $U_{\Stwo}(X)$, over which $\pi$ is an isomorphism, so $\pi$ matches them up with those of $X$. 
\epf


\section{Ubiquity of Cohen--Macaulay blowing ups} \lab{ubiquity}

In general, it is difficult to determine whether the blowing up $\Bl_I(R)$ of an ideal $I \subset R$ in a Noetherian ring is Cohen--Macaulay. For instance, even if $R$ is Cohen--Macaulay, the same need not be true for its blowing up at a maximal ideal (see \cite{HIO88}*{14.11}) in spite of such centers being ``nice,'' in particular, normally flat, equimultiple, normally Cohen--Macaulay, etc. Strikingly, Kawasaki constructed a broad class of ideals $I$ for which $\Bl_I(R)$ is Cohen--Macaulay. His ideals are the backbone of the results of this paper, so the goal of this section is to review them in \Cref{Kaw-blowup}. Such $I$ are built out of the following ``CM-secant'' sequences (see also Remark \ref{p-std-rem}).

\bd \lab{CM-sec-def}
For a finite module $M$ over a Noetherian, local ring $(R, \fm)$, a sequence $r_1, \dotsc, r_s \in \fm$
\benumr
\item \lab{CSD-i}
is \emph{secant for $M$} (in the sense of \cite{BouAC}*{VIII, \S3, no.~2, Def.~1}) if for every $1 \le i \le s$ we have
\[
\qq \dim(\Supp(M/(r_1, \dotsc, r_{i})M)) < \dim(\Supp(M/(r_1, \dotsc, r_{i - 1})M))
\]
(since $\dim(\emptyset) = -\infty$ (see \S\ref{conv}), only the empty sequence is secant when $M = 0$);

\item \lab{CSD-ii}
is \emph{CM-secant for $M$} if it is secant for $M$, the $R$-module $M/(r_1, \dotsc, r_s)M$ is Cohen--Macaulay, and for every $1 \le i \le s$ we have
\[
\qqq r_i \in \prod_{j < \dim(\Supp(M/(r_1,\, \dotsc,\, r_{i - 1})M)) } \Ann_R\p{H^j_\fm(R, M/(r_1, \dotsc, r_{i - 1})M)}
\]
(so that, informally and imprecisely, $r_i$ vanishes to a high enough order on the non-Cohen--Macaulay locus of $M/(r_1, \dotsc, r_{i - 1})M$; see \Cref{Ann-ref}~\ref{AR-b} for a precise such statement).
\eenum
\ed

\brems
\remi \lab{sec-def}
A sequence is secant for $M$ if and only if it is a part of a sequence of parameters for $M$: indeed, for every $1 \le i \le s$, we have
\[
\qqq \dim(\Supp(M/(r_1, \dotsc, r_{i - 1})M)) - 1 \le \dim(\Supp(M/(r_1, \dotsc, r_{i})M)),
\]
so $r_1, \dotsc, r_s$ is secant if and only if 
\[
\qq \dim(\Supp(M/(r_1, \dotsc, r_{s})M)) = \dim(\Supp(M)) - s.
\]
Thus, any permutation of a secant for $M$ sequence $r_1, \dotsc, r_s$ is still secant for $M$, as is any $r_1, \dotsc, r_{i - 1}, r_i r_j$ with $1 \le i \le  j \le s$ (neither $r_i$ nor $r_j$ vanishes at any point of $\Supp(M/(r_1, \dotsc, r_{i - 1})M)$ of maximal coheight).  Any $M$-regular sequence is~secant if $M \neq 0$.

\remi \lab{Sch-ann}
By a result of Schenzel \cite{Sche82}*{2.4.2}, the product ideals in the definition of a CM-secant sequence control the failure of the sequence to be regular as follows: for a finite module $M$ over a Noetherian, local ring $(R, \fm)$ and any secant for $M$ sequence $r_1, \dotsc, r_s \in \fm$ with $s \ge 1$, the $r_s$-torsion of $M/(r_1, \dotsc, r_{s - 1})M$ is killed by the ideal 
\[
\qq \prod_{j < \dim(\Supp(M))} \Ann_R(H^j_\fm(R, M))
\]
(if $M$ is Cohen--Macaulay, then this reconfirms  that any such $r_1, \dotsc, r_s$ is $M$-regular).

\remi \lab{p-std-rem}
If $r_1, \dotsc, r_s$ is a system of parameters for $M$, then the $0$-dimensional $M/(r_1, \dotsc, r_s)M$ is always Cohen--Macaulay. Thus, in this case, $r_1, \dotsc, r_s$ is CM-secant for $M$ if and only if $r_s, \dotsc, r_1$ is a $p$-standard system of parameters for $M$ in the sense of \cite{NTC95}*{2.4} or, equivalently, a $p$-standard system of parameters of type $s - 1$ in the sense of \cite{Kaw00}*{2.6}. Conversely, any CM-secant for $M$ sequence $r_1, \dotsc, r_s$ can be extended to a CM-secant sequence of parameters for $M$ because $M/(r_1, \dotsc, r_s)M$ is, by assumption, Cohen--Macaulay.
\erems

Our next goal is the aforementioned \Cref{Ann-ref} that will be key for the constructions in \S\ref{local-case}.

\bpp[Dualizing complexes] \lab{dual-comp}
We recall from \cite{SP}*{\href{https://stacks.math.columbia.edu/tag/0A7B}{0A7B}} that for a Noetherian ring $R$, a complex $\omega_R^\bullet$ of $R$-modules is \emph{dualizing} if its cohomology modules are finitely generated, its image in the derived category $D(R)$ is isomorphic to a finite complex of injective $R$-modules (so $\omega_R^\bullet \in D^b_\mathrm{coh}(R)$), and 
\[
R \isomto R\Hom_R(\omega_R^\bullet, \omega_R^\bullet) \q \x{in} \q D(R).
\]
The resulting property that justifies the name `dualizing' is \cite{SP}*{\href{https://stacks.math.columbia.edu/tag/0A7C}{0A7C}}: namely, the functor
\[
R\Hom_R(-, \omega_R^\bullet) \colon D^b_\mathrm{coh}(R) \ra D^b_\mathrm{coh}(R) \q \x{is an involutive antiequivalence of categories.}
\]
For us, $\omega_R^\bullet$ will only be important through its image in $D(R)$, so we will abuse terminology and call the later a \emph{dualizing complex}, even without choosing an actual complex representing it.

If $R$ has a dualizing complex $\omega_R^\bullet$, then the latter is unique in $D(R)$ up to tensoring with an object that Zariski locally on $R$ is a shift of $R$ (see \cite{SP}*{\href{https://stacks.math.columbia.edu/tag/0A7F}{0A7F}}). The formation of $\omega_R^\bullet$ commutes with localization (see \cite{SP}*{\href{https://stacks.math.columbia.edu/tag/0A7G}{0A7G}}), and also with quotients as follows: for any $R \surjects R'$, the object 
\[
R\Hom_R(R', \omega_R^\bullet) \in D(R')
\]
is dualizing (see 
 \cite{SP}*{\href{https://stacks.math.columbia.edu/tag/0A7I}{0A7I}}).  Thus, if $R$ is local with residue field $k$, then 
\[
R\Hom_R(k, \omega_R^\bullet) \q \x{is isomorphic to $k$ placed in some degree $n \in \bZ$,}
\]
and we say that $\omega_R^\bullet$ is \emph{normalized} if this $n$ is $0$ (so, in general, $\omega_R^\bullet[n]$ is normalized).

By \cite{Kaw02}*{1.4} (which settled Sharp's conjecture \cite{Sha79}*{4.4}), a Noetherian ring $R$ has a dualizing complex if and only if it is a quotient of a finite dimensional Gorenstein ring. We will only use the simpler `if' implication \cite{SP}*{\href{https://stacks.math.columbia.edu/tag/0DW7}{0DW7}}, in particular, that any complete, Noetherian, local $R$ has a dualizing complex. As follows from the above criterion and Remark \ref{quasi-CM} (and is better seen directly \cite{SP}*{\href{https://stacks.math.columbia.edu/tag/0AWY}{0AWY}, \href{https://stacks.math.columbia.edu/tag/0DW9}{0DW9}, and \href{https://stacks.math.columbia.edu/tag/0A80}{0A80}}), an $R$ that has a dualizing complex is CM-excellent. If, in addition, $R$ is local, then it is Cohen--Macaulay if and only if $\omega_R^\bullet$ is concentrated in a single degree \cite{SP}*{\href{https://stacks.math.columbia.edu/tag/0AWS}{0AWS}}.
\epp

\blem \lab{Ann-ref}
Let $(R, \fm)$ be a Noetherian local ring that has a normalized dualizing complex $\omega^\bullet_R$.
\benum
\item \lab{AR-a}
For every finite $R$-module $M$ and every $j \in \bZ$,
\[
\qq \Ann_R(H^j_\fm(R, M)) = \Ann_R(H^{-j}(R\Hom_R(M, \omega_R^\bullet))).
\]

\item \lab{AR-b}
For every finite $R$-module $M$ with equidimensional support, the ideal
\[
\qq \prod_{j < \dim(\Supp(M))} \Ann_R(H^j_\fm(R, M)) = \prod_{j < \dim(\Supp(M))}  \Ann_R(H^{-j}(R\Hom_R(M, \omega_R^\bullet)))
\]
cuts out a closed subscheme of $\Spec(R)$ whose complement is the Cohen--Macaulay locus of~$M$.
\eenum
\elem

\bpf \hfill
\benum
\item
Let $E$ be the injective hull over $R$ of the residue field $R/\fm$. The local duality theorem \cite{Har66}*{V,~6.2} (or \cite{Sche82}*{1.3.4}) gives an identification
\[
\qq R\Gamma_\fm(R, M) \cong R\Hom_R(R\Hom_R(M, \omega^\bullet_R), E),
\]
which is natural in $M$. Since $E$ is injective, the functor $\Hom_R(-, E)$ is exact and the previous identification gives
\[
\qq H^j_\fm(R, M) \cong \Hom_R( H^{-j}(R\Hom_R(M, \omega_R^\bullet)), E).
\]
Thus, it suffices to show that for every finite $R$-module $M'$, the inclusion
 \[
\qq  \Ann_R(M') \subset \Ann_R(\Hom_R(M', E))
 \] 
 is an equality. By the Nakayama lemma, a nonzero $M'$ always has $R/\fm$ as a quotient, and $R/\fm \subset E$. Thus, $\Hom_R(M', E) = 0$ is equivalent to $M' = 0$. Since the functor $\Hom_R(-, E)$ is exact, it remains to apply it to each sequence 
\[
\qq 0 \ra M'\langle r \rangle \ra M' \xra{r} M' \ra M'/rM' \ra 0.
\]

\item
For a prime $\fp \subset R$ of coheight $\delta(\fp)$, by \cite{SP}*{\href{https://stacks.math.columbia.edu/tag/0A7Z}{0A7Z}}, the dualizing for $R_\fp$ complex 
\[
\qq (\omega^\bullet_R|_{R_\fp})[-\delta(\fp)]
\]
is normalized. Moreover, due to the finite generation of the $R$-modules 
\[
\qq H^{-j}(R\Hom_R(M, \omega_R^\bullet)),
\]
the formation of the product of their annihilator ideals commutes with localization at $\fp$. Thus, a $\fp \in \Supp(M)$ does not belong to the closed subset in question if and only if
\[
\qq H^{-j}(R\Hom_{R_\fp}(M_\fp, (\omega^\bullet_R|_{R_\fp})[-\delta(\fp)])) = 0
\]
for every 
\[
\qq j < \dim(\Supp(M)) - \delta(\fp) = \dim(\Supp(M_\fp)),
\]
where we used the equidimensionality of $\Supp(M)$ and the catenarity of $R$ (see \S\ref{dual-comp}) for the last equality. By \cite{SP}*{\href{https://stacks.math.columbia.edu/tag/0A7U}{0A7U}}, this vanishing amounts to the Cohen--Macaulayness of $M_\fp$.
\qedhere
\eenum
\epf

\brem
More generally, if the support of $M$ in \Cref{Ann-ref}~\ref{AR-b} is arbitrary, then, by \cite{Sche82}*{2.4.6}, the product ideal in question cuts out a closed subscheme that set-theoretically equals
\[
\qq \{ \fp \in \Supp(M)\, \vert\, \depth_{R_\fp}(M_\fp) + \dim(R/\fp) < \dim(\Supp(M)) \} \subset \Spec(R)
\]
and, in particular, contains each irreducible component of $\Supp(M)$ of nonmaximal dimension. 
\erem

Kawasaki established many pleasant properties of CM-secant sequences, some of which will be reviewed in \Cref{Kaw-input}. In particular, he proved that they satisfy the following weak version, which originates with Huneke \cite{Hun82}, of the definition of a regular sequence.

\bd \lab{d-seq-def}
For a module $M$ over a commutative ring $R$, a sequence $r_1, \dotsc, r_s \in R$ is a \emph{$d$-sequence for $M$} if for every $i \ge 1$ scaling by $r_i$ is injective on the submodule 
\[
\tst \f{(r_1,\, \dotsc,\, r_s) M}{(r_1,\, \dotsc,\, r_{i - 1})M} \subset \f{M}{(r_1,\, \dotsc,\, r_{i - 1})M}.
\] 
\ed

\brems
\remi \lab{d-seq-other-def}
The definition implies that the $r_i$-torsion and the $r_{i}r_{i'}$-torsion submodules of 
\[
\qq \tst \f{M}{(r_1,\, \dotsc,\, r_{i - 1})M}
\]
agree for every $1 \le i \le i' \le s$. This property is often taken as the definition of a $d$-sequence. The resulting notion agrees with that of \Cref{d-seq-def}, see \cite{HIO88}*{38.6 b) and its proof}. 

\remi \lab{Goto-Shimoda}
Set $\fr \ce (r_1, \dotsc, r_s) \subset R$. \Cref{d-seq-def} says that for $i \ge 1$ no element of $\fr M$ maps to a nonzero $r_i$-torsion element of $\f{M}{(r_1,\, \dotsc,\, r_{i - 1})M}$. This gives the first of the equalities
\be \lab{Goto-Shimoda-eq}
\qqq ((r_1, \dotsc, r_{i - 1})M :_M r_i) \cap \fr^nM = (r_1, \dotsc, r_{i - 1})M \cap \fr^nM = (r_1, \dotsc, r_{i - 1})\fr^{n - 1}M
\ee
that hold for every $n \ge 1$ and $1 \le i \le s$: indeed, by increasing induction on $n$ and decreasing induction on $i$ the second one follows  from its trivial cases $i = s + 1$ and $n = 1$ as follows. Only ``$\subset$'' needs an argument and for $n \ge 2$ and $i \le s$ the inductive assumptions give
\[
\qq  (r_1, \dotsc, r_{i - 1})M \cap r_i \fr^{n - 1}M \overset{\ref{d-seq-def}}{=} r_i ( (r_1, \dotsc, r_{i - 1})M \cap \fr^{n - 1}M) = r_i(r_1, \dotsc, r_{i - 1})\fr^{n - 2}M.
\]
Consequently, they also give both inclusions in the desired 
\[
\qq (r_1, \dotsc, r_{i - 1})M \cap \fr^nM \subset (r_1, \dotsc, r_{i - 1})M \cap (r_1, \dotsc, r_{i})\fr^{n - 1}M \subset (r_1, \dotsc, r_{i - 1})\fr^{n - 1}M.
\]
The equality \eqref{Goto-Shimoda-eq} is due to Goto and Yamagishi and is also proved in \cite{Kaw00}*{2.2}.
\erems

The ``amplification by induction'' technique of the proof of the equality \eqref{Goto-Shimoda-eq}, whose $n = 1$ case amounts to a definition, is emblematic of the arguments that go into the following key result. 

\bprop[Kawasaki] \lab{Kaw-input}
For a finite module $M$ over a Noetherian local ring $(R, \fm)$, any initial subsequence $r_1, \dotsc, r_s$ of a CM-secant for $M$ sequence $r_1, \dotsc, r_s, r_{s + 1}, \dotsc, r_{\wt{s}}$ with $s \le \wt{s}$ satisfies the following conditions \ref{KI-b}--\ref{KI-a}, where  we set 
\[
\tst I_j \ce \prod_{i = 1}^j (r_1, \dotsc, r_i) \subset R.
\]
\benumr

\item \lab{KI-b}
For any secant for $M/(r_1, \dotsc, r_s)M$ sequence $r'_1, \dotsc, r'_{s'} \in \fm$ with $M' \ce M/(r_1', \dotsc, r'_{s'})M$, 
\[
\qq r_s, \dotsc, r_1 \q \x{is a $d$-sequence for $M'$.}
\]

\item \lab{KI-c}
For any secant for $M/(r_1, \dotsc, r_s)M$ sequence $r'_1, \dotsc, r'_{s'} \in \fm$ with $M' \ce M/(r_1', \dotsc, r'_{s'})M$,  
\[
\qqq \x{ for $m, n > 0$, the $r_s^m$-torsion and the $(r_1, \dotsc, r_s)$-torsion of} \q M'/I_{s - 1}^nM' \q \x{coincide}
\]
\up{when $s \le 1$, the claim holds vacuously}, in particular, this $r_s^m$-torsion does not depend on $m$.

\item \lab{KI-a}
For any secant for $M/(r_1, \dotsc, r_s)M$ sequence $r'_1, \dotsc, r'_{s'} \in \fm$ with $M' \ce M/(r'_1, \dotsc, r'_{s' - 1})M$, if $r_{s'}$ is $(M'/(r_1, \dotsc, r_s)M')$-regular, then it is $M'$-regular and $(M'/I_s^nM')$-regular for $n > 0$.

\eenum
\eprop

\bpf
Assertion \ref{KI-b} is a special case of \cite{Kaw00}*{2.10} (see also Remark \ref{p-std-rem}) and is the shortest one to prove: one deduces it from a slightly more general \cite{Kaw00}*{2.9}, which is a consequence of the result of Schenzel that we reviewed in Remark \ref{Sch-ann}. In turn, the assertions \ref{KI-c} and \ref{KI-a} are special cases of \cite{Kaw00}*{3.2} and  \cite{Kaw00}*{3.3}, respectively, and lie deeper, even though their proofs do not use any external inputs other than the already mentioned result of Schenzel. These proofs rest on the key \cite{Kaw00}*{3.1}, which presents five statements $(A_{ij}), \dotsc, (E_{ij})$ whose flavor is similar to that of \eqref{Goto-Shimoda-eq}, and then proves them by a somewhat lengthy interwoven induction on the difference $j - i$. The base case $i = j$ comes from the already mentioned \cite{Kaw00}*{2.9--2.10} and the Goto--Yamagishi result \eqref{Goto-Shimoda-eq}. We saw the main technique of the proofs in Remark \ref{Goto-Shimoda}: by an inductive assumption, one first gets a slightly weaker statement and then bootstraps using the definition of a $d$-sequence. 

Even though we do not reproduce the cited proofs here, we stress that they are written clearly and are not difficult to follow: one only needs to read \cite{Kaw00}*{\S2 and \S3} (and we already covered much of \cite{Kaw00}*{\S2}). In addition, we are always in the simplest case of ``$p$-standard sequences of type $d - 1$'' (see Remark \ref{p-std-rem}), so we have no need for ``$d^+$-sequences'' that are relevant for the more general types of $p$-standard sequences.  To aid the reading, we mention some harmless misprints:
\benuma
\item
 in the last line of the proof of \cite{Kaw00}*{2.8}, the `$n_{s + 1}, \dotsc, x_d$'  should be `$n_{s + 1}, \dotsc, n_d$';
 
\item \lab{change-2}
in the line before \cite{Kaw00}*{(3.1.5)}, the `$y_u \in \fa(M/\fq_kM)$' should be `$y_u \in \fa(M/\fq_kM)$ or $y_u \in \fa(M)$'---this is needed in order to be able to apply $(E_{i + 1,\, j})$ at the end of Step 6 of the proof of \cite{Kaw00}*{3.1} and is not necessary when the $p$-standard sequence is of type $d - 1$;\footnote{The same correction is made in \cite{Kaw02}*{3.6}, which contains a slightly more general version of \cite{Kaw00}*{3.1}. }

\item
the change \ref{change-2} should also be made before the last display of Step 3 of the proof of \cite{Kaw00}*{3.1}; 

\item
in \cite{Kaw00}*{(3.1.6)}, the `$y_1, \dotsc, y_{u - 1}$' should be `$y_1, \dotsc, y_u$';
 
\item
in the last line of \cite{Kaw00}*{3.3}, the `$k \le i \le j$' should be `$k \le i \le j \le d$.' \qedhere
\eenum
\epf

\brem \lab{robust}
For a finite module $M$ over a Noetherian local ring $(R, \fm)$, the conditions \ref{KI-b}--\ref{KI-a} are such that if a sequence $r_1, \dotsc, r_s \in \fm$ satisfies them, then it continues to do so once $M$ is replaced by $M/(r_1', \dotsc, r'_{s'})M$ for any secant for $M/(r_1, \dotsc, r_s)M$ sequence $r_1', \dotsc, r'_{s'} \in \fm$.
\erem

\bpp[Blowing up modules] \lab{str-tran}
For a scheme $X$, a quasi-coherent $\sO_X$-module $\sM$, and a quasi-coherent ideal $\sI \subset \sO_X$, we consider the quasi-coherent $\sO_{\Bl_{\sI}(X)}$-module 
\[
\tst \Bl_{\sI}(\sM) \q \x{associated to the graded $(\bigoplus_{n \ge 0} \sI^n)$-module} \q \bigoplus_{n \ge 0} \sI^n\sM.
\]
Concretely, for an affine open $\Spec(R) \subset X$ with 
\[
I \ce \Gamma(R, \sI) \q\x{and}\q M \ce \Gamma(R, \sM),
\]
and an $i \in I$, the homogeneous localization $R_{(i)}$ is the $R$-subalgebra of $R[\f{1}{i}]$ generated by the $\f{i'}{i}$ with $i' \in I$, and $\Gamma(R_{(i)}, \Bl_{\sI}(\sM))$  is the $R_{(i)}$-submodule $M_{(i)} \subset M[\f{1}{i}]$ generated by the image of $M$. In particular, $i$ is $M_{(i)}$-regular, so $\Bl_{\sI}(\sM)$ has no nonzero sections supported on the exceptional divisor.  There is a natural surjection 
\[
\sM|_{\Bl_{\sI}(X)} \surjects \Bl_{\sI}(\sM) 
\]
that corresponds to the multiplication map
\[
R_{(i)}  \tensor_R M  \surjects M_{(i)},
\]
which is an isomorphism away from the vanishing locus of $\sI$. Due to the preceding paragraph, its kernel consists of the local sections of $\sM|_{\Bl_{\sI}(X)}$ supported on the exceptional divisor. Consequently, $\Bl_\sI(\sM)$ is nothing else but the strict transform of $\sM$ in the sense of \cite{RG71}*{I.5.1.1 (ii)}.

By \cite{SP}*{\href{https://stacks.math.columbia.edu/tag/080A}{080A}}, for any additional quasi-coherent ideal $\sI' \subset \sO_X$, we have canonical $X$-isomorphisms
\be \lab{blow-product}
\Bl_{\sI\cdot\sI'}(X) \cong \Bl_{\sI\cdot \sO_{\Bl_{\sI'}(X)}}(\Bl_{\sI'}(X)) \cong \Bl_{\sI'\cdot \sO_{\Bl_{\sI}(X)}}(\Bl_{\sI}(X)).
\ee
The relation with the strict transform gives the corresponding identifications
\be \lab{blow-more}
\Bl_{\sI\cdot\sI'}(\sM) \cong \Bl_{\sI\cdot \sO_{\Bl_{\sI'}(X)}}(\Bl_{\sI'}(\sM)) \cong \Bl_{\sI'\cdot \sO_{\Bl_{\sI}(X)}}(\Bl_{\sI}(\sM)).
\ee
\epp

The following result of Kawasaki is our eventual source of Cohen--Macaulayness.

\bthm[Kawasaki] \lab{Kaw-blowup}
For a finite module $M$ over a Noetherian, local ring $(R, \fm)$ and a CM-secant for $M$ sequence $r_1, \dotsc, r_{\wt{s}} \in \fm$, the product ideal $I \ce \prod_{i = 1}^{\wt{s}} (r_1, \dotsc, r_i)$ is such that
\[
\Bl_I(M) \q \x{is a Cohen--Macaulay module on} \q \Bl_{I}(R).
\]
\ethm

\emph{Proof.} 
We loosely follow the proof of \cite{Kaw00}*{4.1}. Since $\Bl_I(M)$ injects into its restriction to the complement of the exceptional divisor of $\Bl_I(R)$ (see \S\ref{str-tran}), the support of $\Bl_I(M)$ is  the schematic image in $\Bl_I(R)$ of 
\[
\Supp(M) \setminus (\Supp(M) \cap \Spec(R/I)).
\] 
By the universal property of blowing up (or by \cite{SP}*{\href{https://stacks.math.columbia.edu/tag/080E}{080E}} directly), this schematic image is the blowing up of $\Supp(M)$ at the restriction of $I$. In particular, by \cite{HIO88}*{12.14}, its dimension is $\le \dim(\Supp(M))$. Consequently, since Cohen--Macaulayness is stable under localization (see \cite{EGAIV1}*{0.16.5.10 (i)}), it suffices to show that the depth of $\Bl_I(M)$ at each closed point of its support is $\ge \dim(\Supp(M))$.

For flexibility in subsequent reductions, we will argue this under more general assumptions: instead of requiring that $r_1, \dotsc, r_{\wt{s}} \in \fm$ be CM-secant for $M$, we only require that each of its initial subsequences $r_1, \dotsc, r_s$ with $s \le \wt{s}$ satisfy the properties \ref{KI-b}, \ref{KI-c}, and \ref{KI-a} of \Cref{Kaw-input} and that the quotient 
\[
\ov{M} \ce M/(r_1, \dotsc, r_{\wt{s}})M
\]
be Cohen--Macaulay. Due to \Cref{robust}, this property of the sequence $r_1, \dotsc, r_{\wt{s}}$ persists if one replaces $M$ by $M/r'M$ for some $\ov{M}$-regular $r' \in \fm$. Moreover, \Cref{Kaw-input}~\ref{KI-a} and the snake lemma applied to the short exact sequence
\[
\tst 0 \ra \bigoplus_{n \ge 0} I^n M \ra \bigoplus_{n \ge 0} M \ra \bigoplus_{n \ge 0} M/I^nM \ra 0
\]
show that $r'$ is $(\bigoplus_{n \ge 0} I^n M)$-regular with 
\[
\tst (\bigoplus_{n \ge 0} I^n M)/r'(\bigoplus_{n \ge 0} I^n M) \isomto \bigoplus_{n \ge 0} I^n (M/r'M).
\]
Thus, since $r'$ lies in $\fm$, and hence vanishes at every closed point of $\Bl_I(R)$, we may replace $M$ by $M/r'M$ without losing any generality (see Remark \ref{sec-def}). By iterating this process, we reduce to the case when $M/(r_1, \dotsc, r_{\wt{s}})M$ is $0$-dimensional, to the effect that we seek to show that the depth of $\Bl_I(M)$ at every closed point of its support is $\ge \wt{s}$. For this, we set 
\[
\tst I_s \ce \prod_{i = 1}^s (r_1, \dotsc, r_s) \q\x{for}\q 0 \le s \le \wt{s}
\]
and seek to show by induction on $s$ that the depth of $\Bl_{I_s}(M)$ at each closed point of its support is at least $s$. Any such point lies above $\fm$, so the $r_i$ vanish at it.

The base case $s = 0$ is trivial, and the case $s = 1$ follows from \S\ref{str-tran}, which implies that $r_1$ is not a zero divisor on 
\[
\Bl_{I_1}(M) \cong M/M\langle r_1^\infty\rangle.
\]
Thus, we suppose that $s \ge 2$. By \eqref{blow-product}, we have the~$R$-map
\[
\Bl_{(r_1,\, \dotsc,\, r_s)}(\Bl_{I_{s - 1}}(R)) \cong \Bl_{I_s}(R) \ra \Bl_{I_{s - 1}}(R)
\]
and the induced map on local rings
\[
R'' \la R'
\]
to the local ring $R''$ of $\Bl_{I_s}(R)$ at a closed point of $\Supp(\Bl_{I_s}(M))$ in question from the local ring $R'$ of $\Bl_{I_{s - 1}}(R)$ at the image of this point. We let $\fm'' \subset R''$ and $\fm' \subset R'$ be the maximal ideals and let $M''$ and $M'$ be the corresponding stalks of $\Bl_{I_s}(M)$ and $\Bl_{I_{s - 1}}(M)$. Since $M' \neq 0$ (see, for instance, \cite{EGAI}*{9.5.5}), we may assume by induction that for any $M$ that satisfies the assumptions (we induct on $s$ for all possible $M$ simultaneously) we have
\be \lab{Mpr-to-Mprpr}
 \depth_{R'}(M') \ge s - 1, \q \x{and seek to show that} \q \depth_{R''}(M'') \ge s.
\ee
By \eqref{blow-product} once more, 
\[
\Bl_{I_{s - 1}}(R) \cong \Bl_{(r_1,\, \dotsc,\, r_{s - 1})}(\Bl_{I_{s - 2}}(R)),
\]
so there is an $1 \le i \le s - 1$ such that $r_i$ generates the ideal $(r_1, \dotsc, r_{s - 1})R'$ and is both $R'$-regular and $M'$-regular (see \S\ref{str-tran}). Consequently, 
\[
(r_1, \dotsc, r_{s})R' = (r_a, r_b)R'  \q \x{with} \q \{a, b\} = \{i, s\}.
\]
Thus, due to the explicit affine cover of $\Bl_{(r_a, r_b)}(R')$ given by affine blowing up algebras (see \cite{SP}*{\href{https://stacks.math.columbia.edu/tag/0804}{0804}}), we may assume that $R''$ is the localization of the subring $R'[\f{r_a}{r_b}] \subset R'[\f{1}{r_b}]$ at a maximal ideal whose $R'$-preimage is $\fm'$. Explicitly, since each maximal ideal of $(R'/\fm')[T]$ is generated by a monic polynomial, there is a monic $f(T) \in R'[T]$ such that 
\[
R'' \simeq (R'[T]/(r_bT - r_a))_{(\fm',\, f(T))}/(r_b\x{-torsion})
\]
and
\[
 M'' \overset{\eqref{blow-more}}{\simeq} (M' \tensor_{R'} R'')/(r_b\x{-torsion})
\]
with the preimage of $\fm''$ in $R'[T]$ generated by $\fm'$ and $f(T)$. The spectral sequence for local cohomology \cite{SP}*{\href{https://stacks.math.columbia.edu/tag/0BJC}{0BJC}} (with \cite{SP}*{\href{https://stacks.math.columbia.edu/tag/0955}{0955}}), the modified \v{C}ech complex interpretation of local cohomology \cite{SP}*{\href{https://stacks.math.columbia.edu/tag/0A6R}{0A6R}}, and the $r_b$-torsion freeness of $M''$ give the identifications
\be \lab{punch-1}
H^j_{\fm''}(R'', M'') \cong H^{j - 1}_{\fm''}(R'', H^1_{(r_b)}(R'', M'')) \overset{\x{\cite{SP}*{\href{https://stacks.math.columbia.edu/tag/0BJB}{0BJB}}}}{\cong} H^{j - 1}_{(\fm',\, f(T))}(R'[T], H^1_{(r_b)}(R'[T], M'')).
\ee
The \v{C}ech complex interpretation also shows that $H^{\ge 1}_{(r_b)}$ vanishes on $r_b^\infty$-torsion modules, so 
\be \lab{punch-2}
H^1_{(r_b)}(R'[T], M'') \cong H^1_{(r_b)}(R'[T], M' \tensor_{R'} (R'[T]/(r_bT - r_a))_{(\fm',\, f(T))})
\ee
In the last identification, we may take the localization $(-)_{(\fm',\, f(T))}$ out of the cohomology (see \cite{SP}*{\href{https://stacks.math.columbia.edu/tag/0ALZ}{0ALZ}}) and, since the coefficients is an $R'[T]/(r_bT - r_a)$-module and $(r_b) = (r_a, r_b)$ in $R'[T]/(r_bT - r_a)$, replace $H^1_{(r_b)}$ by $H^1_{(r_a,\, r_b)}$. Therefore, the combination of \eqref{punch-1} and \eqref{punch-2} gives the identification
\be\lab{punch-3}
H^j_{\fm''}(R'', M'') \cong H^{j - 1}_{(\fm',\, f(T))}(R'[T], H^1_{(r_a,\, r_b)}(R'[T], (M' \tensor_{R'} R'[T])/(r_bT - r_a))).
\ee
We seek the vanishing of the left side of \eqref{punch-3} for $j < s$ (see \eqref{Mpr-to-Mprpr} and \cite{SGA2new}*{III,~3.3~(iv)}). The right side is in terms of $M'$, so we will deduce the vanishing from the following claims.

\bcl \lab{claim-2}
Both $r_i$ and $r_s$ kill $H^1_{(r_i,\, r_s)}(R', M')$.
\ecl

\bcl \lab{claim-1}
We have $H^j_{\fm'}(R', H^{j'}_{(r_i,\, r_s)}(R', M')) = 0$ for $j < s - 2$ and any $j'$. 
\ecl

To deduce the promised vanishing we will use the sequence
\be \lab{punch-seq}
0 \ra M' \tensor_{R'} R'[T] \xra{r_bT - r_a} M' \tensor_{R'} R'[T] \ra (M' \tensor_{R'} R'[T])/(r_bT - r_a) \ra 0,
\ee
which is exact because either $r_a$ or $r_b$ is $M'$-regular (the one with the index $i$).\footnote{In fact, both $r_a$ and $r_b$ are $M'$-regular, see the proof of \Cref{claim-2}.} By flat base change,
\be \lab{punch-4}
H^j_{\fm'}(R'[T], H^{j'}_{(r_a,\, r_b)}(R'[T], M'\tensor_{R'} R'[T] )) \cong H^j_{\fm'}(R', H^{j'}_{(r_a,\, r_b)}(R', M')) \tensor_{R'} R'[T],
\ee
so scaling by the monic $f(T)$ is injective on these groups. Thus, the spectral sequence \cite{SP}*{\href{https://stacks.math.columbia.edu/tag/0BJC}{0BJC}}~gives
\[
H^{j}_{(\fm',\, f(T))}(R'[T], \wt{M}) \cong H^1_{(f(T))}(R'[T], H^{j - 1}_{\fm'}(R'[T], \wt{M})) \q \x{with} \q \wt{M} = H^{j'}_{(r_a,\, r_b)}(R'[T], M'\tensor_{R'} R'[T]).
\]
By combining this with \eqref{punch-4}, we therefore conclude from \Cref{claim-1} that
\be \lab{punch-5}
H^{j}_{(\fm',\, f(T))}(R'[T], H^{j'}_{(r_a,\, r_b)}(R'[T], M'\tensor_{R'} R'[T])) = 0 \q \x{for} \q j < s - 1 \q \x{and any} \q j'.
\ee
On the other hand, \Cref{claim-2} and the version of \eqref{punch-4} for $H^{j'}_{(r_a,\, r_b)}$ alone imply that 
\be \lab{punch-6}
\x{multiplication by} \q r_aT - r_b \q \x{is the zero map on} \q H^{1}_{(r_a,\, r_b)}(R'[T], M'\tensor_{R'} R'[T] ).
\ee
Moreover, since $(r_a, r_b) = (r_bT - r_a, r_b)$ in $R'[T]$, the \v{C}ech complex interpretation shows the vanishing
\be\lab{punch-7}
H^2_{(r_a,\, r_b)}(R'[T], (M' \tensor_{R'} R'[T])/(r_bT - r_a)) = 0.
\ee
The cohomology with supports $H^*_{(r_a,\, r_b)}(R'[T], -)$ sequence that arises from \eqref{punch-seq} combines with \eqref{punch-5}, \eqref{punch-6}, and \eqref{punch-7} to show the desired vanishing of the left side of \eqref{punch-3} for $j < s$: namely, for $ j < s - 1$ we have
\[
H^{j}_{(\fm',\, f(T))}(R'[T], H^1_{(r_a,\, r_b)}(R'[T], (M' \tensor_{R'} R'[T])/(r_bT - r_a))) = 0.
\]
To conclude the proof we will argue \Cref{claim-2,claim-1} by using the inductive assumption and the properties \ref{KI-b}--\ref{KI-c}, which each initial subsequence of $r_1, \dotsc, r_s$ was assumed satisfy.

\emph{Proof of Claim \uref{claim-2}.} 
The coherent module $\Bl_{I_{s - 1}}(M)$ is associated to the graded $(\bigoplus_{n \ge 0} I_{s - 1}^n R)$-module $(\bigoplus_{n > 0} I_{s - 1}^n M)$: indeed, we may omit the summand in degree $0$ (or in any initial segment of degrees) by \cite{EGAII}*{2.5.4, 2.5.6}. The property \ref{KI-b} ensures that $r_s, \dotsc, r_1$ is a $d$-sequence for $M$, so $r_s$ is $(\bigoplus_{n > 0} I_{s - 1}^n M)$-regular, and hence also $M'$-regular. Thus, by the spectral sequence \cite{SP}*{\href{https://stacks.math.columbia.edu/tag/0BJC}{0BJC}},
\be \lab{bam-2}
 H^1_{(r_i,\, r_s)}(R', M') \cong H^0_{(r_i)}(R', H^1_{(r_s)}(R', M')) \overset{\x{\cite{SP}*{\href{https://stacks.math.columbia.edu/tag/0A6R}{0A6R}}}}{\cong} \varinjlim_{m > 0} H^0_{(r_i)}(R', M'/r_s^mM'),
\ee
where the transition maps are induced by multiplication by $r_s$. We will conclude by arguing that both $r_i$ and $r_s$ kill each $H^0_{(r_i)}(R', M'/r_s^mM')$. For this, we begin with the short exact sequence
\be \lab{bam-1}
0 \ra \bigoplus_{n > 0} \f{(I_{s - 1}^nM) :_M r_s^m}{I_{s - 1}^nM + (0 :_M r_s^m)} \xra{r_s^m} \bigoplus_{n > 0} \f{I_{s - 1}^nM}{r_s^m I_{s - 1}^nM} \ra \bigoplus_{n > 0} \p{ I_{s - 1}^n \p{\f{M}{r_s^m M}}} \ra 0
\ee
of graded $(\bigoplus_{n \ge 0} I_{s - 1}^n R)$-modules. The property \ref{KI-b} ensures that $r_i, \dotsc, r_1$ is a $d$-sequence for $M/r_s^mM$ and $(r_1, \dotsc, r_i)$ is a factor of $I_{s - 1}$, so $r_i$ is $\p{\bigoplus_{n > 0} \p{ I_{s - 1}^n \p{\f{M}{r_s^m M}}}}$-regular (actually, we will only use that $r_i$ is a nonzerodivisor on the stalk at $\fm'$ of the corresponding sheaf, which follows from \S\ref{str-tran}). Consequently, it suffices to show that both $r_i$ and $r_s$ kill 
\[
 \f{(I_{s - 1}^nM) :_M r_s^m}{I_{s - 1}^nM + (0 :_M r_s^m)} \q \x{for}\q n > 0.
\]
This follows from the property \ref{KI-c}: indeed, it implies that scaling by any element of $(r_1, \dotsc, r_s) \subset R$ brings $(I_{s - 1}^nM) :_M r_s^m$ inside $I_{s - 1}^n M$.
\QED

\emph{Proof of Claim \uref{claim-1}.} 
The claim is that the $E_2^{jj'}$-entries with $j < s - 2$ of the spectral sequence 
\[
E_2^{jj'} = H^j_{\fm'}(R', H^{j'}_{(r_i,\, r_s)}(R', M')) \Ra H^{j + j'}_{\fm'}(R', M')
\]
vanish. Due to the \v{C}ech complex interpretation of $H^{j'}_{(r_i,\, r_s)}(R', -)$, this vanishing holds for $j' > 2$. Due to the $M'$-regularity of $r_i$, the same holds for $j' = 0$. Thus, we only need to handle $j' = 1$ and $j' = 2$. By the inductive assumption \eqref{Mpr-to-Mprpr} on $M'$, the abutment $H^{j + j'}_{\fm'}(R', M')$ vanishes for $j + j' < s - 1$, so the case $j' = 2$ alone would suffice. To handle it, similarly to \eqref{bam-2}, we have 
\be \lab{bam-3}
 H^2_{(r_i,\, r_s)}(R', M') \cong H^1_{(r_i)}(R', H^1_{(r_s)}(R', M')) \cong \varinjlim_{m > 0}\p{ H^1_{(r_i)}(R', M'/r_s^mM')}.
\ee
For each $m > 0$, let $M'_{(m)}$ be the stalk at $\fm'$ of the sheaf that corresponds to the graded $(\bigoplus_{n \ge 0} I_{s - 1}^n R)$-module $\bigoplus_{n > 0} \p{ I_{s - 1}^n \p{\f{M}{r_s^m M}}}$. As we saw in the proof of \Cref{claim-2}, due to \ref{KI-c}, scaling by $r_i$ kills the kernel term of the short exact sequence \eqref{bam-1}. Thus, this sequence gives the identification
\be \lab{bam-4}
H^1_{(r_i)}(R', M'/r_s^mM') \cong H^1_{(r_i)}(R', M'_{(m)}) \q \x{for every} \q m > 0.
\ee
By combining \eqref{bam-3} and \eqref{bam-4}, we get
\be\ba \lab{bam-bam}
 H^j_{\fm'}(R', H^2_{(r_i,\, r_s)}(R', M')) &\cong  \varinjlim_{m > 0} H^j_{\fm'}(R', H^1_{(r_i)}(R', M'_{(m)})) \\ &\cong  \varinjlim_{m > 0}\p{\varinjlim_{m' > 0} H^j_{\fm'}(R', M'_{(m)}/r_i^{m'}M'_{(m)})}.
\ea\ee
The inductive assumption \eqref{Mpr-to-Mprpr} applies to each $M'_{(m)}$ (see \Cref{robust}), so 
\[
H^j_{\fm'}(R', M'_{(m)}) = 0 \q \x{for}\q j < s - 1 \q \x{and}\q m > 0.
\]
Consequently, since $r_i$ is $M'_{(m)}$-regular (see \S\ref{str-tran}), we also have 
\[
H^j_{\fm'}(R', M'_{(m)}/r_i^{m'}M'_{(m)}) = 0 \q \x{for}\q j < s - 2 \q\x{and}\q m' > 0,
\]
and the claim follows from \eqref{bam-bam}.
\QEDD



\section{Macaulayfication in the local case} \lab{local-case}

The key step of the inductive Macaulayfication procedure is to treat an $X$ that is projective over a complete Noetherian local ring $(R, \fm)$ with $X \setminus X_{R/\fm}$ already Cohen--Macaulay and $X_{R/\fm} \subset X$ a divisor. In such a situation, \Cref{Mac-local} exhibits a Macaulayfying blowing up $\wt{X} \ra X$ whose center does not meet $X \setminus X_{R/\fm}$. The crucial trick is to regard $X_{R/\fm^n} \subset X$ for $n > 0$ as a hypersurface and to use it in building a center to be blown up whose localizations are as in \Cref{Kaw-blowup}. Since the latter is applied locally, the local principality of $\sI_{X_{R/\fm^n}} \cong (\sI_{X_{R/\fm}})^n$ is enough to make this~work.


\bpp[Biequidimensionality] \lab{biequi}
We recall from \cite{EGAIV1}*{0.14.3.3}\footnote{Recall from \cite{ILO14}*{XIV, \S2.4, footnote (i)} that \cite{EGAIV1}*{0.14.3.3 b)} should read ``$X$ est cat\'{e}naire, \'{e}quidimensionnel, et ses composantes irr\'{e}ductibles sont \'{e}quicodimensionnelles'' (see also \cite{Hei17}*{\S3}).} that a scheme $X$ is \emph{biequidimensional} if it is Noetherian, of finite dimension, and the saturated chains of specializations of its points all have the same length, necessarily equal to $\dim(X)$. Evidently, every biequidimensional $X$ is locally equidimensional (see \S\ref{formal-eq}).  For an example, by \cite{EGAIV2}*{5.2.1}, every connected, locally equidimensional $X$ that is of finite type over a field is biequidimensional. 

If $X$ is biequidimensional, then so is every nowhere dense closed subscheme $X' \subset X$ whose coherent sheaf of ideals is locally principal. Conversely, if such an $X'$ in a Noetherian, catenary, locally equidimensional $X$ contains all the closed points of $X$ and is biequidimensional, then $X$ is also biequidimensional because all of its closed points have the same height equal to $\dim(X') + 1$.
\epp

\bpp[Dualizing complexes on schemes] \lab{DC-schemes}
We recall from \cite{SP}*{\href{https://stacks.math.columbia.edu/tag/0A87}{0A87}} that for a Noetherian scheme $X$, an object
\[
\omega_X^\bullet \in D^b_{\mathrm{coh}}(\sO_X)
\]
is a \emph{dualizing complex} if its restriction to every affine open $\Spec(R) \subset X$ is the image of a dualizing complex in $D(R)$ (see \S\ref{dual-comp}). If $X$ has a dualizing complex, then, by \cite{SP}*{\href{https://stacks.math.columbia.edu/tag/0AA3}{0AA3}} and the Nagata compactification \cite{Del10}*{1.6}, so does every finite type, separated $X$-scheme~$X'$. 

A pleasant case is when a Noetherian $X$ that has a dualizing complex $\omega_X^\bullet$ is  biequidimensional. Then, due to \cite{SP}*{\href{https://stacks.math.columbia.edu/tag/0AWF}{0AWF}}, there is a unique $n \in \Gamma(X, \un{\bZ})$ such that 
for any $x \in X$ with coheight $\delta(x)$, the dualizing for $\sO_{X,\, x}$ complex $(\omega_X^\bullet|_{\sO_{X,\, x}})[-\delta(x) + n(x)]$ is normalized in the sense reviewed in \S\ref{dual-comp}. If $n = 0$, as may be arranged by a unique locally constant shift, then we call $\omega^\bullet_X$ \emph{normalized}. 
\epp

In the proof of \Cref{Mac-local}, we will use the following version of the avoidance lemma.

\blem[\cite{GLL15}*{5.1}] \lab{GLL-input}
For a Noetherian ring $R$, a projective $R$-scheme $X$, a closed subscheme $Y \subset X$ that does not contain any positive dimensional irreducible component of any $R$-fiber of $X$, points $x_1, \dotsc, x_n \in X \setminus Y$, and a very  $R$-ample line bundle $\sL$ on $X$, there are an $N > 0$ and an $f \in \Gamma(X, \sL^{\tensor N})$ such that the closed subscheme of $X$ cut out by the coherent ideal
\[
\tst \im\p{ (\sL^{\tensor N})\i \xra{f} \sO_X } \subset \sO_X
\]
contains $Y$ but does not contain any of the $x_1, \dotsc, x_n$.  \QED
\elem

\bprop \lab{Mac-local}
For a complete, Noetherian, local ring $(R, \fm)$, a locally equidimensional, projective $R$-scheme $X$ such that $X_{R/\fm} \subset X$ is a divisor and $X \setminus X_{R/\fm}$ is Cohen--Macaulay, and finitely many coherent $\sO_X$-modules $\sM$ with $\sM|_{X\setminus X_{R/\fm}}$ Cohen--Macaulay and $\abs{\Supp(\sM)} = \abs{X}$, there is a closed subscheme $Z \subset X_{R/\fm^n}$ for some $n > 0$ such that $\Bl_Z(X)$ is Cohen--Macaulay and its coherent modules $\Bl_Z(\sM)$ \up{see \uS\uref{str-tran}} are also all Cohen--Macaulay. 
\eprop

\bpf
By including $\sO_X$ among $\sM$, we reduce to arranging the Cohen--Macaulayness of~the~$\Bl_Z(\sM)$.

By \cite{EGAIII1}*{5.5.1}, the $R$-properness of $X$ ensures that the $(R/\fm)$-fiber of each connected component of $X$ is connected. By passing to such a component, we assume that both $X$ and $X_{R/\fm}$ are connected (so nonempty). By \Cref{fed-stab}~\ref{FS-a} (with \S\ref{formal-eq}), the divisor $X_{R/\fm}$ inherits local equidimensionality from $X$. Then \S\ref{biequi} ensures that $X_{R/\fm}$, and then also $X$, is biequidimensional.

By \S\ref{DC-schemes} and \S\ref{dual-comp}, the biequidimensional, proper $R$-scheme $X$ has a normalized dualizing complex $\omega_X^\bullet$.
Since 
\[
R\sH om_{\sO_X}(\sM, \omega^\bullet_X)  \in D^b_{\mathrm{coh}}(\sO_X),
\]
for every $\sM$ the following ideal is coherent:
\[
\tst \sI_\sM \ce \prod_{j < \dim(X)} \Ann_{\sO_X}(H^{-j}(R\sH om_{\sO_X}(\sM, \omega^\bullet_X))) \subset  \sO_X.
\]
Due to biequidimensionality, for an $x \in X$ with coheight $\delta(x)$, we have 
\[
\dim(\sO_{X, x}) + \delta(x) = \dim(X),
\]
so
\[
\tst (\sI_\sM)_x 
\overset{\x{\ref{Ann-ref}~\ref{AR-a}}}{=} \prod_{j < \dim(\sO_{X,\,  x})} \Ann_{\sO_{X,\, x}}(H^j_{\fm(x)}(\sO_{X,\,x}, \sM_{x})),
\]
where 
$\fm(x) \subset \sO_{X,\, x}$ is the maximal ideal. In particular, by \Cref{Ann-ref}~\ref{AR-b}, the ideal $\sI_\sM$ cuts out a closed subscheme of $X$ whose open complement is $\CM(\sM)$. Thus, there is an $N_0 > 0$ such that
\[
\sJ_1 \ce (\sI_{X_{R/\fm}})^{N_0} \subset \sI_\sM \subset \sO_X \q \x{for every} \q \sM.
\]
We will construct a decreasing sequence $X \supset X_1 \supset X_2 \supset \dotsc \supset X_d$ of biequidimensional closed subschemes of $X$ of strictly decreasing dimension starting with $X_1 \ce \un{\Spec}(\sO_X/\sJ_1)$ as follows. 

Assume that $X_m$ has already been constructed. If each $\sM|_{X_m}$ is Cohen--Macaulay, then set $d \ce m$ and stop. Otherwise, consider the dualizing complex 
\[
\omega_{X_m}^\bullet \ce R\sH om_{\sO_X}(\sO_{X_m}, \omega_X^\bullet) \q\x{on}\q X_m
\]
that is automatically normalized (see \cite{SP}*{\href{https://stacks.math.columbia.edu/tag/0AX1}{0AX1}}). Analogously to the case of $\sI_\sM$, for every $\sM$ the ideal
\[
\tst \sI_{\sM,\, m} \ce \prod_{j < \dim(X_m)} \Ann_{\sO_{X_m}}(H^{-j}(R\sH om_{\sO_{X_m}}(\sM|_{X_m}, \omega^\bullet_{X_m}))) \subset  \sO_{X_m}
\]
 is coherent and, by \Cref{Ann-ref}~\ref{AR-a}, for every $x \in X_m$, its stalk is
\[
\tst (\sI_{\sM,\, m})_x = 
\prod_{j < \dim(\sO_{X_m,\,  x})} \Ann_{\sO_{X_m,\, x}}(H^j_{\fm(x)}(\sO_{X_m,\,x}, (\sM|_{X_m})_x)).
\]
By \cite{SP}*{\href{https://stacks.math.columbia.edu/tag/0BJB}{0BJB}} (with \cite{SP}*{\href{https://stacks.math.columbia.edu/tag/0955}{0955}}), if we identify $\sI_{\sM,\, m}$ with its preimage in $\sO_{X}$, then this stalk~becomes
\be \lab{new-stalk}
\tst (\sI_{\sM,\, m})_x = 
\prod_{j < \dim(\sO_{X_m,\,  x})} \Ann_{\sO_{X,\, x}}(H^j_{\fm(x)}(\sO_{X,\,x}, (\sM|_{X_m})_x)).
\ee
In particular, by \Cref{Ann-ref}~\ref{AR-b}, the ideal $\sI_{\sM,\, m}$ cuts out a closed subscheme $Y_{\sM,\, m} \subset X_m$ whose complement is $\CM(\sM|_{X_m})$, so that $Y_{\sM,\, m}$ does not contain any generic point of $X_m$. Thus, by \Cref{GLL-input}, for a fixed very $R$-ample line bundle $\sL$ on $X$, there are an $N_m > 0$ and an $f_m \in \Gamma(X, \sL^{\tensor N_m})$ such that the closed subscheme of $X$ cut out by the locally principal ideal
\[
\tst \sJ_{m + 1 } \ce \im\p{(\sL^{\tensor N_m})\i \xra{f_m} \sO_{X}} \subset \sO_X
\]
contains each $Y_{\sM,\, m}$ but does not contain any generic point of $X_m$. We let $X_{m + 1}$ be the scheme-theoretic intersection of $X_m$ and the closed subscheme cut out by $\sJ_{m + 1}$. Since $X_{m + 1}$ is nowhere dense in $X_m$, it inherits biequidimensionality by \S\ref{biequi}. For dimension reasons, the process stops.

For each $x \in X$ and $1 \le m \le d$, let $j_{m,\, x} \in \sO_{X,\, x}$ be a generator of the ideal $(\sJ_{m})_x \subset \sO_{X,\, x}$ and, for the sake of convenience, set $j_{d + 1,\, x} \ce 1 \in \sO_{X,\, x}$. Let $m_x$ be the smallest $0 \le m \le d$ such that $j_{m + 1,\, x} \in \sO_{X,\, x}^\times$. By construction, the sequence $j_{1,\, x}, \dotsc, j_{m_x,\, x}$ is secant for $\sO_{X,\, x}$ (see \Cref{CM-sec-def}~\ref{CSD-i}) and $\sO_{X,\, x}/(j_{1,\, x}, \dotsc, j_{m_x,\, x})$ is the local ring of $X_{m_x}$ at the point $x$. Since $j_{m_x + 1,\, x} \in \sO_{X,\, x}^\times$, we have 
\[
x \in \CM(\sM|_{X_{m_x }}) \q\x{for every}\q \sM.
\]
Consequently, by construction and \eqref{new-stalk}, the sequence
\[
j_{1,\, x}, \dotsc, j_{m_x ,\, x} \in \fm(x) \q \x{is CM-secant for every} \q \sM_{x}
\]
(see \Cref{CM-sec-def}~\ref{CSD-ii}). Thus, since the coherent ideal
\[
\tst \sJ \ce \prod_{i = 1}^{d}(\sJ_1, \dotsc, \sJ_i)
\]
satisfies
\[
\tst \sJ_x = \prod_{i = 1}^{m_x }(j_{1,\, x}, \dotsc, j_{i,\, x}) \subset \sO_{X,\, x} \q \x{for every} \q x \in X,
\]
\Cref{Kaw-blowup} implies that  $\Bl_\sJ(\sM)$ is Cohen--Macaulay for every $\sM$. It remains to set 
\[
Z \ce \un{\Spec}(\sO_X/\sJ) \q\x{and} \q n \ce N_0 \cdot d. \qedhere
\]
\epf


\section{Macaulayfication in the global case} \lab{global-case}

We are ready to use \Cref{Mac-local} and Noetherian induction to obtain a global Macaulayfication.

\blem\lab{RG-lemma}
For a Noetherian scheme $S$, a closed subscheme $Z \subset S$, the blowing up $\Bl_Z(S) \xra{b} S$, and a closed subscheme $Z' \subset \Bl_Z(S)$, there is a closed subscheme $Z'' \subset S$ such that 
\[
\abs{Z''} = \abs{Z} \cup b(\abs{Z'})
\]
and there are $S$-identifications
\[
\Bl_{Z''}(S) \cong \Bl_{Z'}(\Bl_Z(S))
\]
and, for any quasi-coherent $\sO_S$-module $\sM$
\[
 \Bl_{Z''}(\sM) \cong \Bl_{Z'}(\Bl_Z(\sM)).
\]
\elem

\bpf
The claim is a special case of \cite{SP}*{\href{https://stacks.math.columbia.edu/tag/080B}{080B}}, except for the aspect about the modules. The latter follows by interpreting both $\Bl_{Z''}(\sM)$ and $\Bl_{Z'}(\Bl_Z(\sM))$ as strict transforms, see \S\ref{str-tran}.
\epf

\bprop \lab{Mac-main}
For every CM-excellent, locally equidimensional, Noetherian scheme $X$ and finitely many coherent $\sO_X$-modules $\sM$ with $\abs{\Supp(\sM)} = \abs{X}$, there is a closed subscheme 
\[
\tst Z \subset X \q\x{that is disjoint from the dense open \up{see \uS\uref{open-Sn}}}\q \CM(X) \cap \p{\bigcap_\sM \CM(\sM)}
\]
such that $\Bl_{Z}(X)$ is Cohen--Macaulay, and its coherent modules $\Bl_{Z}(\sM)$ are also all Cohen--Macaulay.
\eprop

\bpf
The claim is evidently true in the case when $X$ and $\sM$ are all Cohen--Macaulay: one chooses $Z = \emptyset$. In general, given an open subscheme $U \subsetneq X$ that contains 
\[
\tst \CM(X) \cap \p{\bigcap_\sM \CM(\sM)}
\]
and for which the claim holds with some closed subscheme $Z \subset U$, we need to argue that the claim also holds for some strictly larger open. For this, we first use \cite{EGAI}*{9.4.7} to extend $Z$ to a closed subscheme $\wt{Z} \subset X$. We will show that there is a closed subscheme of $\Bl_{\wt{Z}}(X)$ that is disjoint from 
\[
(\Bl_{\wt{Z}}(X))|_U \cong \Bl_Z(U),
\]
whose blowing up is Cohen--Macaulay over an open neighborhood of a fixed generic point of $X \setminus U$, and for which the corresponding strict transforms of the modules $\Bl_{\wt{Z}}(\sM)$ are also all Cohen--Macaulay over this neighborhood. By \Cref{RG-lemma}, this will allow us to enlarge $U$.

For the remaining claim about $\Bl_{\wt{Z}}(X)$, due to \cite{EGAI}*{9.4.7} again and a limit argument based on the openness of the Cohen--Macaulay loci of CM-excellent schemes (see \S\ref{open-Sn}), we may localize at the fixed generic point of $X \setminus U$ to assume that $X$ is local and $U$ is the complement of the closed point. Then $X = \Spec(R)$ for a Noetherian, local $(R, \fm)$ and to simplify notation, we let $Y$ be the projective $R$-scheme $\Bl_{\wt{Z}}(X)$ that comes equipped with the finitely many coherent modules 
\[
\wt{\sM} \ce \Bl_{\wt{Z}}(\sM)
\]
of set-theoretically maximal support (see \S\ref{str-tran}). By assumption, the open subscheme $Y_U$ is Cohen--Macaulay, and so are its coherent modules $\wt{\sM}|_{U}$. We seek a closed subscheme $Z \subset Y$ that is disjoint from $Y_U$ such that $\Bl_Z(Y)$ and $\Bl_Z(\wt{\sM})$ are all Cohen--Macaulay, so, thanks to \Cref{RG-lemma}, we may precompose $Y \ra \Spec(R)$ with $\Bl_{Y_{R/\fm}}(Y) \ra Y$ to assume that 
\[
Y_{R/\fm} \subset Y \q\x{is a divisor.}
\]
By assumption, $X$ is formally equidimensional (see \S\ref{formal-eq}), so, by \Cref{fed-stab}~\ref{FS-b}, the blowing up $Y_{\wh{R}}$ of $\wh{R}$ is locally equidimensional. 
By \Cref{Sn-descend}, the base change $Y_{U_{\wh{R}}}$ and its modules $\wt{\sM}|_{Y_{U_{\wh{R}}}}$ inherit Cohen--Macaulayness. Thus, \Cref{Mac-local} applies to 
\[
Y_{\wh{R}} \ra \Spec(\wh{R})
\]
and gives a closed subscheme $Z \subset Y_{R/\fm^n}$ for some $n > 0$ such that $\Bl_{Z}(Y_{\wh{R}})$ and its coherent modules $\Bl_Z(\wt{\sM}|_{Y_{\wh{R}}})$ are all Cohen--Macaulay. Since blowing up commutes with flat base change, 
\[
\Bl_{Z}(Y_{\wh{R}}) \cong (\Bl_Z(Y))_{\wh{R}} \q \x{and} \q \Bl_Z(\wt{\sM}|_{Y_{\wh{R}}}) \cong \Bl_Z(\wt{\sM})|_{Y_{\wh{R}}},
\]
so \Cref{Sn-descend}  implies that $\Bl_Z(Y)$ and its modules $\Bl_Z(\wt{\sM})$ are all Cohen--Macaulay, as desired. 
\epf

\bthm \lab{main-thm-pf}
For every CM-quasi-excellent, Noetherian scheme $X$ equipped with finitely many coherent $\sO_X$-modules $\sM$ with $\abs{\Supp(\sM)} = \abs{X}$, there are a composition
\[
\wt{X} \ce \Bl_{Z}(X') \ra X' \xra{\pi'} X
\]
and for each $\sM$ a coherent $\sO_{X'}$-module $\sM'$ for which $\abs{\Supp(\sM')} = \abs{X'}$ such that $\wt{X}$ is Cohen--Macaulay, its coherent modules $\Bl_{Z}(\sM')$ are also all Cohen--Macaulay, and
\benumr
\item
$X'$ is CM-excellent and locally equidimensional\uscolon

\item
$\pi'$ is finite, birational, and is an isomorphism over the open 
\[
\qq \tst U \ce U_{\x{\Stwon}}(X) \cap \p{\bigcap_\sM U_{\x{\Stwon}}(\sM)} \subset X
\]
\up{see \uS\uref{open-Sn}} that is dense in both $X$ and $X'$ and for which $\sM'|_U \cong \sM|_U$\uscolon

\item
$Z \subset X'$ is a closed subscheme that is disjoint from the dense open 
\[
\qq\tst U' \ce \CM(X') \cap \p{\bigcap_\sM \CM(\sM')};
\]
\item \label{4}
$U'$ is also dense in $\Bl_Z(X')$, so that, in particular, the map $\wt{X} \ra X$ is birational\uscolon
\eenum
 if $X$ itself is CM-excellent and locally equidimensional, then we may choose $X' = X$ and $\sM' = \sM$.
\ethm

\bpf 
\Cref{S2-ify} supplies an \Stwon, finite, birational $X$-scheme $X'$ with $X'|_U \isomto U$ such that $U$ is dense both in $X$ and $X'$, as well as \Stwon, coherent $\sO_{X'}$-modules $\sM'$ such that 
\[
\sM'|_U \cong \sM|_U \q\x{and}\q \abs{\Supp(\sM')} = \abs{X'}.
\]
By \Cref{miracle-catenary} (with \S\ref{formal-eq}), the CM-quasi-excellent, \Stwo scheme $X'$ is CM-excellent and locally equidimensional. If $X$ itself is CM-excellent and locally equidimensional, then instead we choose 
\[
X' \ce X \q\x{and}\q \sM' \ce \sM.
\] 
\Cref{Mac-main} applies to $X'$ and supplies a closed subscheme $Z \subset X'$ that is disjoint from $U'$
such that $\Bl_{Z}(X')$ and its coherent modules $\Bl_{Z}(\sM')$ are all Cohen--Macaulay. The generic points of $\Bl_Z(X')$ lie away from the divisor given by the preimage of $Z$, so $U'$ is also dense in $\Bl_Z(X')$. 
\epf

\brem \lab{rem-proj}
The map $\wt{X} \ra X$ in \Cref{main-thm-pf} is projective, see \cite{EGAII}*{5.5.5~(ii) and 6.1.11}. 
\erem

We conclude with the following converse result, an analogue of Grothendieck's \cite{EGAIV2}*{7.9.5}.

\bprop \lab{converse}
Let $X$ be a locally Noetherian scheme. If for every integral, closed subscheme $X' \subset X$ there are a Cohen--Macaulay scheme $\wt{X'}$ and a proper map $\wt{X'} \ra X'$ that is an isomorphism over a nonempty open subscheme of $X'$ \up{so that $\wt{X'}$ is a Macaulayfication of $X'$}, then $X$ is CM-quasi-excellent.
\eprop

\bpf
The argument is similar to that of \emph{loc.~cit.} Namely, the assumption implies that each $X'$ contains a nonempty open that is Cohen--Macaulay, so we only need to check that the local rings of $X$ have Cohen--Macaulay formal fibers. For this, we may assume that $X = \Spec(R)$ for a local ring $R$ and then, since the formal fibers of a local $X$ are exhausted by those of its irreducible components, that $X$ is also integral. By passing to the closure of the point of $X$ at which the formal fiber is taken, we then reduce further to checking that the generic formal fiber of $X$ is Cohen--Macaulay. 

Let $\pi \colon \wt{X} \ra X$ be a proper map with $\wt{X}$ Cohen--Macaulay and $\pi$ an isomorphism over a nonempty open subscheme of $X$. By the latter condition, it suffices to show that $\wt{X}_{\wh{R}}$ is Cohen--Macaulay. Thus, since $\wt{X}_{\wh{R}}$ inherits CM-excellence, and so also the openness of the Cohen--Macaulay locus (see \S\ref{open-Sn}), from $\wh{R}$, it suffices to show that the local rings of $\wt{X}_{\wh{R}}$ at the closed points  are Cohen--Macaulay. Due to the $X$-properness of $\wt{X}$, these closed points lie over the closed point of $\Spec(\wh{R})$, and analogously for $\wt{X}$, so they are identified with those of $\wt{X}$. Moreover,  for a closed point $\wh{x} \in \wt{X}_{\wh{R}}$  with the image $x \in \wt{X}$, by \cite{EGAIV2}*{7.9.3.1}, the completions $\wh{\sO}_{\wt{X}_{\wh{R}},\, \wh{x}}$ and $\wh{\sO}_{\wt{X},\, x}$ are identified. Since Cohen--Macaulayness ascends to and descends from the completion (see \cite{EGAIV2}*{6.3.5}), the claim follows.
\epf


\end{document}